\documentclass[lettersize,journal]{IEEEtran}
\usepackage{amsmath,amsfonts}

\usepackage{array}
\usepackage[justification=centering]{caption}
\usepackage[caption=false,font=normalsize,labelfont=sf,textfont=sf]{subfig}

\usepackage{textcomp}
\usepackage{stfloats}
\usepackage{url}
\usepackage{verbatim}
\usepackage{graphicx}
\usepackage{cite}
\usepackage{multirow}
\usepackage{longtable}

\usepackage{amsmath}
\usepackage{amsthm}
\usepackage{amssymb}
\usepackage{mathtools}
\usepackage{xspace}

\usepackage{float}
\usepackage{color}
\newcommand{\red}[1]{#1}   % 黒字で表示

\usepackage{enumitem}
\usepackage{fancyhdr}
\usepackage{hyperref}

\usepackage{bm}
\usepackage{algorithm}
\usepackage{algorithmicx}
\usepackage[noend]{algpseudocode}
\algnewcommand\LeftComment[2]{%
\hspace{#1\algindent}$\triangleright$ \eqparbox{COMMENT}{#2} \hfill 
}
\algnewcommand\algorithmicto{\textbf{to}}
\algnewcommand\algorithmicin{\textbf{in}}
\algnewcommand\algorithmicpara{\textbf{parallel}}
\algdef{SE}[FOR]{ForTo}{EndForTo}[2]
    {\algorithmicfor\ #1\ \algorithmicto\ #2\ \algorithmicdo}
    {\algorithmicend\ \algorithmicfor}
\algdef{SE}[FOR]{ForIn}{EndForIn}[2]
    {\algorithmicfor\ #1\ \algorithmicin\ #2\ \algorithmicdo}
    {}
\algdef{SE}[FOR]{ForToPara}{EndForToPara}[2]
    {\algorithmicfor\ #1\ \algorithmicto\ #2\ \algorithmicdo\ \algorithmicin\ \algorithmicpara}
    {\algorithmicend\ \algorithmicfor}
\algtext*{EndForTo}
\algtext*{EndForIn}
\algtext*{EndForToPara}

\makeatletter
\renewenvironment{proof}[1][\proofname]{
  \vskip 2\topsep 
  \hrule
  \vskip0.4\topsep 
  \normalfont  #1\@addpunct{.}\hskip \labelsep\ignorespaces 
  \vskip0.4\topsep 
  \hrule 
  \normalfont 
}{
  \pushQED{\qed} 
  \popQED 
  \vskip\topsep 
  \hrule 
  \vskip\topsep
}
\makeatother

\usepackage[most]{tcolorbox}
\usepackage{lipsum} % ダミーテキスト

\tcbset{
  commentbox/.style={
    colframe=blue!50!black,   % 線の色
    colback=blue!10,          % 背景色
    left=2mm,                 % 左側余白
    boxrule=1pt,              % 線の太さ
    sharp corners,
    before skip=5pt,
    after skip=5pt,
  }
}

\begin{document}

\title{Annealed Mean Field Descent Is Highly Effective for Quadratic Unconstrained Binary Optimization}

% \author{Anonymous}
\author{Kyo~Kuroki,
Thiem~Van~Chu, \IEEEmembership{Member, IEEE},
Masato~Motomura, \IEEEmembership{Fellow, IEEE},
and
Kazushi~Kawamura, \IEEEmembership{Member, IEEE}

\thanks{This work was supported in part by JSPS KAKENHI, Japan, under Grant JP23H05489; in part by JST PRESTO, Japan, under Grant JPMJPR23P1; and in part by JST BOOST, Japan, under Grant JPMJBS2430.}
\thanks{Kyo~Kuroki and Masato~Motomura are with the Institute of Integrated Research, Institute of Science Tokyo, Yokohama 226-8501, Japan.}
\thanks{Thiem~Van~Chu is with the School of Engineering, Institute of Science Tokyo, 152-8550, Japan.}
\thanks{Kazushi~Kawamura is with the Faculty of Science and Engineering, Waseda University, Tokyo 169-8555, Japan.}
}

\maketitle

\begin{abstract}
In recent years, formulating various combinatorial optimization problems as Quadratic Unconstrained Binary Optimization (QUBO) has gained significant attention as a promising approach for efficiently obtaining optimal or near-optimal solutions. While QUBO offers a general-purpose framework, existing solvers often struggle with performance variability across different problems.
This paper (i) theoretically analyzes Mean Field Annealing (MFA) and its variants, which are representative QUBO solvers, and shows that their underlying self-consistent equations do not necessarily correspond to the minimization condition of the Kullback--Leibler divergence between the mean-field approximation and the exact distribution. Moreover, the paper (ii) proposes a novel method, Annealed Mean Field Descent (AMFD), which directly minimizes this divergence to address this limitation.
Through extensive experiments on five benchmark combinatorial optimization problems (Maximum Cut Problem, Maximum Independent Set Problem, Traveling Salesman Problem, Quadratic Assignment Problem, and Graph Coloring Problem), we demonstrate that AMFD exhibits superior performance in many cases and reduced problem dependence compared to state-of-the-art QUBO solvers and Gurobi---a state-of-the-art versatile mathematical optimization solver not limited to QUBO.
\end{abstract}


\section{Introduction}\label{sec:introduction}
\subsection{Background}
\begin{figure}[t]
    \centering
    \includegraphics[width=1\linewidth]{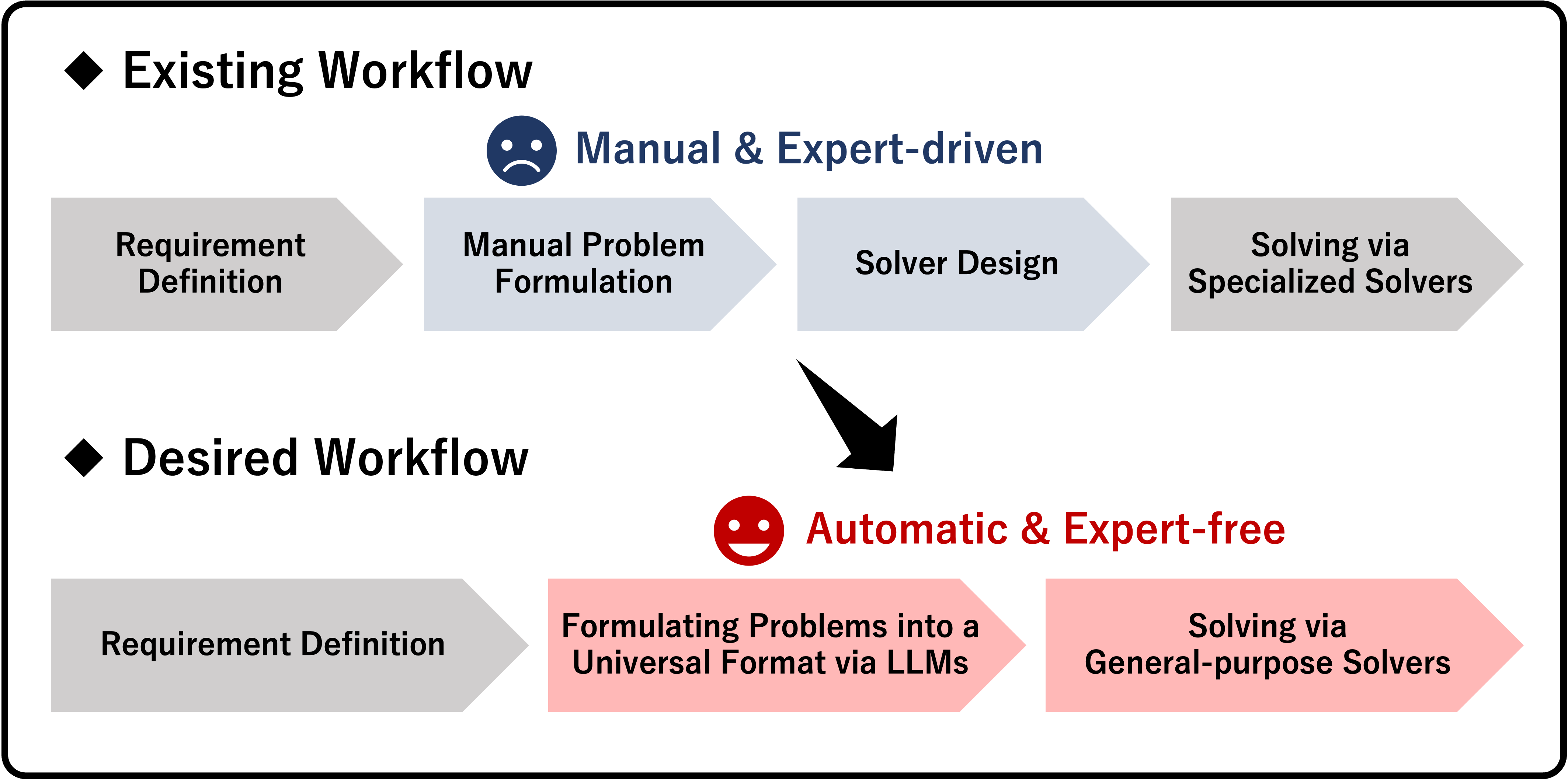}
    \caption{An approach to democratize combinatorial optimization by eliminating formulation and solver design.}
    \label{fig:automation}
\end{figure}

\IEEEPARstart{C}{ombinatorial} optimization problems (COPs), defined as the problem of finding the optimal solution among a huge number of alternatives and constraints, have been applied to real-world challenges such as logistics, scheduling, network design, and resource allocation. Since such problems are often difficult to solve in polynomial time (NP-hard), it is essential to develop efficient algorithms. A common approach to COPs is to develop a specialized algorithm for a specific COP. However, real-world problems frequently deviate from the predefined COPs, such as when the number of constraints increases. In such cases, the specialized algorithms may not be directly applicable and have to be rebuilt. Furthermore, when considering undefined new COPs, it is necessary to examine not only the formulation but also the solver design simultaneously. Unfortunately, formulating COPs and designing solvers is a challenging task for most people. Therefore, removing these barriers is essential to democratize combinatorial optimization (\textbf{Fig.~\ref{fig:automation}}). While automatic formulation is beyond the scope of this paper, recent research has attempted to automate the formulation process using Large Language Models (LLM) \cite{AhmadiTeshnizi2024OptiMUSSO_formulation_by_LLM}. In contrast, this work aims to develop a general-purpose combinatorial optimization solver that obviates the need to develop tailor-made solvers for individual problems.

A general-purpose solver is typically designed to handle Mixed Integer Linear Programming (MILP) problems---problems with mixed integer variables where only linear terms exist in the objective function or constraints, or Mixed Integer Quadratic Programming (MIQP) problems—problems with mixed integer variables where quadratic terms exist in the objective function or constraints.
These formats can accommodate the formulation of most COPs.
A well-known efficient exact solution method for MILP and MIQP uses a technique called the branch-and-cut method\cite{Elf2001_branch_and_cut_method}, which combines the branch-and-bound\cite{branch_and_bound_original, branch_and_bound_principles} and cutting-plane\cite{cutting_plane_original, cutting_plane_for_IP} methods.
The branch-and-bound method divides the problem into subproblems by fixing variables, solves their continuous relaxations to find bounds, and prunes subproblems if their bounds are worse than the best current solution.
The cutting-plane method repeatedly solves a continuous relaxation problem and adds new constraints to make the non-integer variables in that solution meet integer conditions.

However, for COPs defined by formulas that are difficult to refine the lower and upper bounds and add effective constraints in the branch-and-cut method, it is very time-consuming to reach a high-quality solution. Thus, Quadratic Unconstrained Binary Optimization (QUBO) solvers have recently gained significant attention as they maintain the general-purpose nature and are expected to provide high-quality solutions more rapidly than exact methods for MILP or MIQP. 
(Note that most COPs that can be expressed in MILP or MIQP formulation---except for problems that do not involve continuous variables---can be converted to QUBO \cite{Lucas14_Formulation, glover2019tutorial}.)
The motivation for converting COPs into the QUBO form lies in its compatibility with a broad class of efficient metaheuristics--such as Simulated Annealing (SA) \cite{Kirkpatrick83_SA} and Genetic Algorithm (GA) \cite{GA_Holland}--as well as learning-based methods \red{\cite{Schuetz2022_PIGNN, Redwan_RGNN, pugacheva2024QRF_GNN, Toenshoff2019RUN_CSP, Chen2024VCM, Pu2024_DRL}}, both of which are often difficult to apply directly to MILP or MIQP formulations. 
% Learning-based methods for solving QUBO problems employ unsupervised training using a loss function derived either from the QUBO objective evaluated on the solutions predicted by the network’s output layer \cite{Schuetz2022_PIGNN, Redwan_RGNN, pugacheva2024QRF_GNN, Toenshoff2019RUN_CSP}, or from the discrepancy between the predicted solution and the current best solution \cite{Chen2024VCM}. 
\red{Learning-based QUBO solvers either directly minimize the QUBO objective as an unsupervised loss on predicted solutions \cite{Schuetz2022_PIGNN, Redwan_RGNN, pugacheva2024QRF_GNN, Toenshoff2019RUN_CSP}, or learn the solution search strategy by exploiting objective value changes as loss signals \cite{Chen2024VCM} or reinforcement learning rewards \cite{Pu2024_DRL}.}
However, such models often incur substantial training cost and struggle to generalize across diverse classes of COPs, which frequently necessitates retraining for each new problem type. Consequently, the overall computational burden of applying these models to newly defined problems can become prohibitively high.
\red{Furthermore, additional considerations in network design are also needed.}
In contrast, metaheuristics provide training-free optimization methods and, in recent years, have become increasingly compatible with a variety of specialized hardware platforms. Examples include FPGA- and ASIC-based accelerators \cite{Goto_21_bSB, Yamamoto20_STATICA, Kawamura23_RPA, Matsubara20_DA_ASIC, Kashimata2023_SB_ASIC, Yamaoka16_CMOS}, as well as quantum computing systems \cite{Johnson_11_QA_DWave} and optical computing architectures \cite{Inagaki_16_CIM, Honjo_21_CIM_100000spin}. Among such metaheuristics, this work focuses on parallel SA, which offers strong performance without any training cost.

% Some of these heuristics can be executed efficiently by leveraging specialized hardware, such as FPGAs~\cite{Goto_21_bSB}, ASICs~\cite{Yamamoto20_STATICA, Kawamura23_RPA, Matsubara20_DA_ASIC, Kashimata2023_SB_ASIC, Yamaoka16_CMOS}, quantum computers~\cite{Johnson_11_QA_DWave}, and optical computers~\cite{Inagaki_16_CIM, Honjo_21_CIM_100000spin}, thereby accelerating the solution process.
% \red{
% In addition to these approaches, learning-based methods for solving QUBO problems have recently shown considerable promise.
% These techniques typically perform unsupervised training using a loss function derived either from the QUBO objective evaluated on the solutions predicted by the network’s output layer \cite{Schuetz2022_PIGNN, Redwan_RGNN, pugacheva2024QRF_GNN, Toenshoff2019RUN_CSP}, or from the discrepancy between the predicted solution and the current best solution \cite{Chen2024VCM}.
% However, such models often require significant computational effort and struggle to generalize across diverse classes of COPs, necessitating retraining for each new problem type.
% Consequently, the computational cost of applying these approaches to newly defined problems can become prohibitively high.
% In contrast to these learning-based approaches, this work focuses on a promising optimization algorithm: parallel SA, which requires no training and offers strong performance.}

\vspace{-0.0em}
\subsection{Previous Works}
SA simulates a physical phenomenon where the ground state (i.e., the lowest energy state) is obtained by gradually cooling from a high-temperature state to a low-temperature state---whereas, SA inspired by quantum annealing\cite{Kadowaki_98_QA} aims the ground state by gradually weakening the transverse magnetic field from a strong state. 
In applying SA to QUBO problems, binary variables represent the spin orientations of a spin glass model \cite{Kirkpatrick1978_Ising_spinglass} on statistical physics, quadratic coefficients represent interactions between spins, linear coefficients represent external magnetic fields, and the evaluation function is interpreted as the energy of the spin glass system.
By simulating this physical process, we aim to obtain a high-quality solution (i.e., low energy state) to the QUBO problem.
SA itself was invented a long time ago but has received renewed attention in recent years due to the development of parallel SA\cite{Okuyama19_MA, Yamamoto20_STATICA, Kawamura23_RPA, Goto_19_SB, Goto_21_bSB, king18_NMFA, Tiunov_19_simCIM, Kuroki_VMFA}.
The classical SA for QUBO problems follows the process of updating one spin variable per iteration (i.e., serial update).
On the other hand, the parallel SA can update all spin variables simultaneously in one iteration.
% as shown in \textbf{Fig.~\ref{fig:parallel_SA}}.
% \begin{figure}[t]
%     \centering
%     \includegraphics[width=0.9\linewidth]{Figures/difference.png}
%     \caption{The difference between serial SA and parallel SA. The left illustrates serial SA, where a single variable is updated sequentially, whereas the right shows parallel SA, where multiple variables are updated simultaneously.}
%     \label{fig:parallel_SA}
% \end{figure}
Consequently, it has been reported that parallel SA has achieved significantly faster convergence to high-quality solutions by using parallel computers, outperforming not only classical SA but also advanced computing systems such as a quantum annealing system and an optical computing system \cite{Goto_21_bSB, Zeng2024_QA_and_QAIA_benchmark}. Furthermore, since parallel SA can utilize semiconductor digital computers, it has also the advantage of being more scalable than systems using quantum computers or optical computers.

Representative parallel SAs include those that extend Markov Chain Monte Carlo (MCMC) based SA \cite{Kirkpatrick83_SA} to parallel updates, such as Momentum Annealing (MA)\cite{Okuyama19_MA}, Stochastic Cellular Automata (SCA)\cite{Yamamoto20_STATICA}, and Ratio-controlled Parallel Annealing (RPA)\cite{Kawamura23_RPA}. Other examples involve extensions of mean-field-approximation-based SA, commonly referred to as Mean Field Annealing (MFA) \cite{Bilbro88_MFA}, to parallel updates, such as Noisy MFA (NMFA) \cite{king18_NMFA} and Verlet MFA (VMFA) \cite{Kuroki_VMFA}. Additionally, there are methods that simulate binary variables by replacing them with continuous variables, such as Simulated Bifurcation (SB)\cite{Goto_19_SB} and its variants--ballistic SB (bSB) and discrete SB (dSB)\cite{Goto_21_bSB}.

\vspace{0.5em}
\subsection{Challenges and Contributions}
Although various parallel SAs have been developed as state-of-the-art (SOTA) QUBO solvers, they still have the challenge that the solution quality is highly dependent on the COP, and it is difficult to obtain a high-quality solution for certain types of COPs. As a result, despite being general-purpose solvers, the range of COPs they can be practically applied to is limited. 
We believe the factors for this difficulty could be as follows: MCMC-based SAs require a vast number of sampling times to approach the desired probability distribution, MFA-based SAs incur errors due to the mean-field approximation, and methods that replace binary variables with continuous variables introduce errors when replacing. 
Among these three types of SAs, we focus on MFA-based approaches, since their performance is expected to improve as the mean-field approximation becomes more accurate.
Then, this paper actually identifies further room for improvement in the mean field approximation accuracy of existing MFA-based SAs and proposes the Annealed Mean Field Descent (AMFD), which is designed to achieve this improvement.
In experiments, we compare AMFD against five COPs---Maximum Cut Problem (MCP), Maximum Independent Set Problem (MISP), Traveling Salesman Problem (TSP), Quadratic Assignment Problem (QAP), and Graph Coloring Problem (GCP)---all of which have well-established benchmark sets, and demonstrate that AMFD is less problem-dependent than SOTA parallel SAs. Additionally, we benchmark AMFD against Gurobi \cite{gurobi}, a leading commercial solver that employs the branch-and-cut method and various heuristics and techniques to provide high-quality solutions quickly for MILP, MIQP, and QUBO. The experimental results show that AMFD outperformed Gurobi in many cases, despite not combining other heuristics or techniques.

\vspace{1em}
\section{Preliminary}
\subsection{Quadratic Unconstrained Binary Optimization (QUBO)}
The evaluation function of a QUBO problem, which corresponds to the energy of the Ising model in physics, can be expressed as 
\begin{align} \label{eq:QUBO_function}
    H(\boldsymbol{s}) &= \sum_{i=1}^{N_{\text{spin}}}h_{i}s_{i} +  \sum_{i=1}^{N_{\text{spin}}}\sum_{j<i}^{N_{\text{spin}}}Q_{i,j}s_{i}s_{j} + C \notag\\
    &= \boldsymbol{h}^\top\boldsymbol{s} + \frac{1}{2} \boldsymbol{s}^\top\boldsymbol{Q}\boldsymbol{s} + C
\end{align}
where $i,j \in \left\{k|k\in \mathbb{N} \land k \le N_\text{spin}\right\}$ represents spin index, $N_{\text{spin}} \in \mathbb{N}$ is the number of spins, $\boldsymbol{s} \in \{0,1\}^{N_\text{spin}}$ represents the spin orientation, $Q_{i,j} \in \mathbb{R}$ denotes the interaction between $s_i$ and $s_j$, $h_i \in \mathbb{R}$ represents the external bias, and $C \in \mathbb{R}$ is a constant.
The interaction matrix $\boldsymbol{Q}$ is symmetric, and its diagonal elements are set to $0$ (i.e., $Q_{i,j} = Q_{j,i}$ for $1\le i, j\le N_{\text{spin}}$, and $Q_{i,i} = 0$ for $1\le i\le N_{\text{spin}}$). 

The QUBO problems are formulated such that a smaller value of $H(\boldsymbol{s})$ indicates a better solution, meaning the spin state $\boldsymbol{s}$ that minimizes $H(\boldsymbol{s})$ corresponds to the optimal solution. For COPs with constraints, these can be converted into QUBO form by adding a penalty term for variable combinations that violate the constraints.

\subsection{Simulated Annealing (SA)}
SA is a typical method for solving QUBO problems. This is an algorithm that describes a phenomenon in which a metal enters a stable state (i.e., a low-energy state) when the temperature is gradually lowered from a high state to a low state. At high temperatures, thermal fluctuations cause the metal to assume various energy states, implying a global search, while at low temperatures, the metal assumes a low energy state, implying a local search.

From the statistical physics theory\cite{gibbs1902_stastical_mechanics}, the probability of taking spin configuration $\boldsymbol{s}$ at a certain temperature $T$ (including Boltzmann's constant), called the canonical distribution $P_{\text{C}}(\boldsymbol{s})$, is expressed by
\begin{equation}\label{eq:canonical_distribution}
    P_{\text{C}}(\boldsymbol{s}) = \frac{1}{Z}\exp\left(-\frac{H(\boldsymbol{s})}{T}\right),
\end{equation}
where $Z$, called the partition function, is a constant that adds up all spin configurations and is expressed by 
\begin{equation}\label{eq:partition_function}
    Z = \sum_{\boldsymbol{s}} \exp\left(-\frac{H(\boldsymbol{s})}{T}\right).
\end{equation}

However, this probability distribution is difficult to calculate directly because the partition function requires summing over all spin configurations, and the number of combinations grows exponentially as the number of variables increases.
Hence, in order to gain this probability distribution, it is common to use Markov Chain Monte Carlo (MCMC) methods in SA, such as the Metropolis algorithm\cite{Hasting57_MCMC} or Gibbs sampling\cite{Geman84_Gibbs_sampling}. 

Although it has been proven that MCMC methods converge to the true distribution after a sufficient number of sampling times, they require long sampling times, making it challenging to obtain high-quality solutions within realistic timeframes. To improve the convergence speed, recent advances have introduced parallel annealing algorithms using MCMC\cite{Okuyama19_MA, Yamamoto20_STATICA, Kawamura23_RPA} that can update all variables simultaneously in contrast to typical SA which can update only one variable at a time through an iterative process. This approach is known to significantly accelerate convergence when implemented on parallel computing system. Nevertheless, experimentally, MCMC-based SA makes it difficult to achieve high-quality solutions within a realistic time in certain cases.

\subsection{Mean Field Annealing}
Another method that does not use MCMC in SA is Mean Field Annealing (MFA)\cite{Bilbro88_MFA}. In MFA, the canonical distribution is approximated by mean field approximation since it cannot be calculated directly. In mean field approximation, since all variables are assumed to be independent, assumed probability distribution $P_\text{MF}\left(\boldsymbol{s}\right)$ can be expressed as the product of the probability for each spin variable, as shown below:
\begin{equation}\label{eq:mean_field_approximation}
    P_\text{MF}\left(\boldsymbol{s}\right) = \prod_{i=1}^{N_\text{spin}}{p_i\left({{s_i}}\right)}, 
\end{equation}
where $p_i\left({{s_i}}\right)$ means the probability that the $i$-th spin variable becomes $s_i$. 
Accordingly, $x_{i} \in \left\{r|r\in \mathbb{R} \land 0 \le r \le 1\right\}$, the expected value of the $i$-th spin in mean field approximation (called mf-spin), is defined as
\begin{align} \label{eq:mean_spin}
    x_{i} &\stackrel{\mathrm{def}}{=} \sum_{s_{i}=0}^{1} s_{i}p_i\left(s_i\right) = 0\cdot p_i(0) + 1\cdot p_i(1)=p_i(1)
\end{align}
Then, the energy function is approximated using the mean fields---meaning the average force applied to a spin by the interaction with other spins and the external magnetic field---as shown below:
\begin{equation}\label{eq:hamiltonian_approximation}
    H\left(\boldsymbol{s}\right) \simeq \sum_{i=1}^{N_{\text{spin}}} \Phi_{i}s_{i},
\end{equation}
where $\Phi_i$ represents the mean field applied to the $i$-th spin and is expressed by
\begin{equation}\label{eq:mean_field}
    \Phi_{i} = h_i + \sum_{j=1}^{N_\text{spin}} Q_{i,j} x_{j}.
\end{equation}
Substituting \textbf{Eq.~\ref{eq:hamiltonian_approximation}} into \textbf{Eq.~\ref{eq:canonical_distribution}} to approximate the canonical distribution,
\begin{align}
    P_\text{C}\left(\boldsymbol{s}\right) &\simeq \frac{\exp{\left(-\frac{\sum_{i=1}^{N_\text{spin}}\Phi_i s_i}{T}\right)}}{\sum_{\boldsymbol{s}}\exp{\left(-\frac{\sum_{i=1}^{N_\text{spin}}\Phi_i s_i}{T}\right)}} \notag\\
    &=\prod_{i=1}^{N_\text{spin}}\frac{\exp{\left(-\frac{\Phi_{i}s_i}{T}\right)}}{\sum_{s_i=0}^{1}\exp\left(-\frac{\Phi_i s_i}{T}\right)} \notag\\
    \label{eq:approximate_distribution}
    &=\prod_{i=1}^{N_\text{spin}}\frac{\exp{\left(-\frac{\Phi_{i}s_i}{T}\right)}}{1+\exp{\left(-\frac{\Phi_{i}}{T}\right)}}
\end{align}
holds. Here, \textbf{Eq.~\ref{eq:approximate_distribution}} is nothing but the assumed probability distribution $P_\text{MF}\left(\boldsymbol{s}\right)$, so
\begin{equation}\label{eq:spin_probability}
    p_i\left(s_i\right) = \frac{\exp{\left(-\frac{\Phi_{i}s_i}{T}\right)}}{1+\exp{\left(-\frac{\Phi_{i}}{T}\right)}}
\end{equation}
holds. Then, substituting \textbf{Eq.~\ref{eq:spin_probability}} into \textbf{Eq.~\ref{eq:mean_spin}}, $x_i$ can be expressed by the following self-consistent equation:
\begin{equation}\label{eq:self_consistent_equation}
    x_i = \frac{1}{1+\exp{\left(\frac{\Phi_{i}}{T}\right)}}.
\end{equation}
MFA is an algorithm that iteratively updates the mf-spin value according to \textbf{Eq.~\ref{eq:self_consistent_equation}} and is shown in \textbf{Alg.~\ref{Alg:MFA}}.
Unlike MCMC, MFA has no guarantee that it necessarily converges to the true distribution, but convergence is often faster than MCMC, and experimental results often provide better quality solutions than MCMC-based SA.

\begin{algorithm}[t]
\caption{Mean Field Annealing (MFA)} \label{Alg:MFA}
\normalsize
\begin{algorithmic}[1]
\Require $\bm{x}(0)$, $\bm{Q}$, $\bm{h}$, $N_{\text{step}}$, $T(t)$
\Ensure  $\bm{x}(N_{\text{step}})$
\ForTo{$t = 1$}{$N_{\text{step}}$}
    \State Random Spin Selection: $i \in \left\{k|k\in \mathbb{N} \land k \le N_\text{spin}\right\}$
    \State Mean Field Calculation: ${\Phi}_i \gets h_{i} + \sum_{j}Q_{i,j}x_{j}(t-1) $
    \State Spin Update $x_i(t) \gets \frac{1}{1+\exp\left({\frac{{\Phi}_i}{T(t)}}\right)}$
\EndForTo
\ForToPara{$i = 1$}{$N_{\text{spin}}$}
    \If {$x_i(N_{\text{step}}) \geq 0.5$}
         \State  $x_i(N_\text{step}) \gets 1$  
    \Else 
        \State $x_i(N_\text{step}) \gets 0$ 
    \EndIf 
\EndForToPara
\end{algorithmic}
\end{algorithm}

As shown in \textbf{Alg.~\ref{Alg:MFA}}, while MFA usually updates one variable at a time, parallel MFA—which can update all variables simultaneously—has been proposed to speed up convergence as with MCMC-based SA. In parallel MFA, the mean fields for all spins are obtained using the matrix operation shown in 
\begin{equation}
    \boldsymbol{\Phi} = \boldsymbol{h} + \boldsymbol{Q}\boldsymbol{x},
\end{equation}
and all variables simultaneously transition to the next state by the weighted average of the current value and the value obtained in \textbf{Eq.~\ref{eq:self_consistent_equation}}:
\begin{equation}\label{eq:parallel_MFA}
    \boldsymbol{x}_\text{new} = (1-\alpha)\boldsymbol{x}_\text{old} + \alpha\frac{1}{1+\exp{\left(\frac{\boldsymbol{\Phi}_\text{old}}{T}\right)}},
\end{equation}
where $\alpha \in \mathbb{R}$ represents the percentage of weighting.
Noisy MFA (NMFA)\cite{king18_NMFA} adds a noise term to \textbf{Eq.~\ref{eq:parallel_MFA}} to make it easier to escape from local solutions, while Verlet MFA (VMFA)\cite{Kuroki_VMFA} reformulate \textbf{Eq.~\ref{eq:parallel_MFA}} as motion equation and simulate using Verlet method, realizing speed up convergence.

\section{Methods}

\subsection{Theoretical Analysis of MFA}

Here, we conduct a theoretical analysis to improve MFA. In MFA, the energy function is approximated by \textbf{Eq.~\ref{eq:hamiltonian_approximation}}, and the canonical distribution is approximated by \textbf{Eq.~\ref{eq:approximate_distribution}}. However, there might be a better approximation method. Therefore, we consider the Kullback-Leibler (KL) divergence, a measure of the difference between the assumed distribution $P_\text{MF}\left(\boldsymbol{s}\right)$ and the canonical distribution $P_\text{C}\left(\boldsymbol{s}\right)$. Then, the problem of finding the best approximate distribution in terms of the KL divergence metric reduces to the following optimization problem:
\begin{align}
    \label{eq:minimization_KL}
    \arg\min_{P_\text{MF}\left(\boldsymbol{s}\right)} & D_\text{KL}\left(P_\text{MF}\left(\boldsymbol{s}\right) \mid P_\text{C}\left(\boldsymbol{s}\right)\right), \\
    \label{eq:constraint_condition}
    \text{s.t. }& \sum_{s_{i}=0}^{1} p_i\left({s_i}\right) = 1.
\end{align}
\textbf{Eq.~\ref{eq:minimization_KL}} means finding $P_\text{MF}\left(\boldsymbol{s}\right)$ such that the KL divergence between the assumed distribution and the true distribution is as small as possible, and
\textbf{Eq.~\ref{eq:constraint_condition}} implies the constraint that all states add up to 1 from the definition of the probability distribution.
Calculating the KL divergence leads to the following equation:
\begin{align}\label{eq:KL_divergence}
    &D_\text{KL}\left(P_\text{MF}\left(\boldsymbol{s}\right) \mid P_\text{C}\left(\boldsymbol{s}\right) \right) \notag\\
    &=\sum_{i=1}^{N_\text{spin}}\sum_{s_i=0}^{1}p_i(s_i)\ln\left(p_i(s_i)\right)
    + \ln\left(Z\right) \notag\\
    &\quad+ \frac{1}{T} \left(\sum_{i=1}^{N_\text{spin}}h_ix_i + \sum_{i=1}^{N_\text{spin}}\sum_{j<i}Q_{i,j}x_ix_j\right),
\end{align}
which is proved in \textbf{Pr.1}.
The first term in \textbf{Eq.~\ref{eq:KL_divergence}} implies mean-field approximated entropy (called mf-entropy), the second term implies constant, and the third term implies mean-field approximated energy (called mf-energy).
The KL divergence indicates that mf-entropy makes a large contribution when the temperature is high, while the mf-energy (i.e., objective function) makes a large contribution when the temperature is low.
By the method of Lagrange multipliers, the candidate solutions to the constrained minimization problem 
of \textbf{Eq.~\ref{eq:minimization_KL}} and \textbf{Eq.~\ref{eq:constraint_condition}} can be obtained from
\begin{equation}\label{eq:lagrange_multiplier}
    \begin{cases}
      L\left(\boldsymbol{p({s})}, \boldsymbol{\lambda}\right)=&D_\text{KL}\left(P_\text{MF}\left(\boldsymbol{s}\right)\mid P_\text{C}\left(\boldsymbol{s}\right)\right) \\
      &- \sum_{i=1}^{N_\text{spin}}\lambda_i
      \left(1 - \sum_{s_{i}=0}^{1} p_i\left({s_i}\right)\right)\\
      \frac{\partial L\left(\boldsymbol{p({s})},\boldsymbol{\lambda}\right)}{\partial p_i(s_i)} = 0, & \forall{i} \in \left\{k|k\in \mathbb{N} \land k \le N_\text{spin}\right\}\\
      \frac{\partial L\left(\boldsymbol{p({s})},\boldsymbol{\lambda}\right)}{\partial \lambda_i} = 0, & \forall{i} \in \left\{k|k\in \mathbb{N} \land k \le N_\text{spin}\right\}
    \end{cases},
\end{equation}
where $\boldsymbol{\lambda} \in \mathbb{R}^{N_\text{spin}}$ represents Lagrange multiplier and differentiation by $p_i(s_i)$ denotes the functional derivative.
If we solve this, 
\begin{align}\label{eq:spin_probability_2}
    &p_i\left(s_i\right) 
    = \frac{\exp{\left(-\frac{\Phi_{i}s_i}{T}\right)}}{1+\exp{\left(-\frac{\Phi_{i}}{T}\right)}}\\
    &\because \frac{\partial x_i}{\partial p_i(s_i)} = \frac{\partial}{\partial p_i(s_i)}\sum_{s_i=0}^{1}s_ip_i(s_i) = s_i \notag
\end{align}
can be obtained, which is consistent with \textbf{Eq.~\ref{eq:spin_probability}}.
Thus, it is clear that $p_i(s_i)$ in \textbf{Eq.~\ref{eq:spin_probability}} implies being under extreme or saddle point conditions (i.e., \textbf{Eq.~\ref{eq:approximate_distribution}} implies the candidate of $\arg\min_{P_\text{MF}\left(\boldsymbol{s}\right)} D_\text{KL}\left(P_\text{MF}\left(\boldsymbol{s}\right) \mid P_\text{C}\left(\boldsymbol{s}\right)\right)$) because the KL divergence is generally not convex with respect to $\boldsymbol{p(s)}$.

% \clearpage
\begin{proof}[\textbf{Proof 1} : Proof of Eq.~\ref{eq:KL_divergence}]\label{pr:KL_mean_field}
    \quad\par
    \quad\par From the definition of KL divergence,
    \small
    {
        \begin{align}
            &D_\text{KL}\left(P_\text{MF}\left(\boldsymbol{s}\right) \mid P_\text{C}\left(\boldsymbol{s}\right)\right) \notag\\
            &= \sum_{\boldsymbol{s}} P_\text{MF}\left(\boldsymbol{s}\right) 
            \ln\left(\frac{P_\text{MF}\left(\boldsymbol{s}\right)}{P_\text{C}\left(\boldsymbol{s}\right)}\right) \notag\\
            &= \sum_{\boldsymbol{s}} \left[ 
                \left(\prod_{i=1}^{N_{\text{spin}}} p_i\left(s_i\right)\right) 
                \ln\left( 
                    \frac{\prod_{i=1}^{N_{\text{spin}}} p_i\left(s_i\right)}
                         {\frac{1}{Z} \exp\left(-\frac{H(\boldsymbol{s})}{T}\right)}
                \right)
            \right] \notag\\
            \label{eq:KL_difinition}
            &= \sum_{\boldsymbol{s}} \left[ 
                \left(\prod_{i=1}^{N_{\text{spin}}} p_i\left(s_i\right)\right) 
                \left(
                    \sum_{i=1}^{N_{\text{spin}}} \ln\left(p_i\left(s_i\right)\right) 
                    + \ln\left(Z\right) 
                    + \frac{1}{T} H\left(\boldsymbol{s}\right)
                \right)
            \right].
        \end{align}
        }
        \normalsize
     \quad\par \quad Here, \par
             \begin{align}
            &\sum_{\boldsymbol{s}} \left[\left(\prod_{i=1}^{N_{\text{spin}}} p_i\left(s_i\right)\right)\sum_{i=1}^{N_{\text{spin}}} \ln\left(p_i\left(s_i\right)\right)
            \right]\notag\\
            &= \sum_{\boldsymbol{s}}            
                \left[
                    \sum_{i=1}^{N_{\text{spin}}} \left(\prod_{j\ne i}p_j\left(s_j\right)\right)p_i\left(s_i\right)\ln\left(p_i\left(s_i\right)\right) 
                \right]\notag\\
            &=\sum_{i=1}^{N_{\text{spin}}} \left(\prod_{j\ne i}\sum_{s_j=0}^{1}p_j\left(s_j\right)\right)\sum_{s_i=0}^{1}p_i\left(s_i\right)\ln\left(p_i\left(s_i\right)\right) \notag\\      \label{eq:mean_field_approximated_entropy}
            &=\sum_{i=1}^{N_{\text{spin}}}\sum_{s_i=0}^{1}p_i\left(s_i\right)\ln\left(p_i\left(s_i\right)\right),
        \end{align}
        \quad\par \par
             \begin{align}
            &\sum_{\boldsymbol{s}} \left(\prod_{i=1}^{N_{\text{spin}}} p_i\left(s_i\right)\right)\ln(Z)\notag\\
            &=\left(\prod_{i=1}^{N_\text{spin}}\sum_{s_i=0}^{1}p_i\left(s_i\right)\right)\ln(Z)\notag\\
            \label{eq:free_energy}
            &= \ln(Z) ,
        \end{align}
        \quad\par and\par
        \begin{align}
            &\sum_{\boldsymbol{s}} \left[\left(\prod_{i=1}^{N_{\text{spin}}} p_i\left(s_i\right)\right)H(\boldsymbol{s})\right]\notag\\
            &=\sum_{\boldsymbol{s}} \left[\left(\prod_{i=1}^{N_{\text{spin}}} p_i\left(s_i\right)\right)\left(\sum_{i=1}^{N_\text{spin}}h_is_i +\sum_{i=1}^{N_\text{spin}}\sum_{j<i}Q_{i,j}s_is_j\right)\right]\notag\\
            &= \sum_{\boldsymbol{s}} \left(\prod_{j\ne i} p_j\left(s_j\right)\sum_{i=1}^{N_\text{spin}}h_is_i p_i\left(s_i\right)\right)\notag\\
            &\quad+\sum_{\boldsymbol{s}} \left(\prod_{k\ne i,j} p_k\left(s_k\right)\sum_{i=1}^{N_\text{spin}}\sum_{j<i}Q_{i,j}s_ip_i\left(s_i\right)s_jp_j\left(s_j\right)\right)\notag\\ 
            &= \left(\prod_{j\ne i} \sum_{s_j=0}^{1}p_j\left(s_j\right)\right)\left(\sum_{i=1}^{N_\text{spin}}h_i\sum_{s_i=0}^{1}s_i p_i\left(s_i\right)\right)\notag\\
            &\quad+\sum_{\boldsymbol{s}}\left(\prod_{k\ne i,j} p_k\left(s_k\right)\sum_{i=1}^{N_\text{spin}}\sum_{j<i}Q_{i,j}s_ip_i\left(s_i\right)s_jp_j\left(s_j\right)\right)\notag\\
            &=\sum_{i=1}^{N_\text{spin}}h_ix_i+\left(\prod_{k\ne i,j} \sum_{s_k=0}^{1}p_k\left(s_k\right)\right)\notag\\        &\quad\cdot\left(\sum_{i=1}^{N_\text{spin}}\sum_{j<i}Q_{i,j}\sum_{s_i=0}^{1}s_ip_i\left(s_i\right)\sum_{s_j=0}^{1}s_jp_j\left(s_j\right)\right)\notag\\
            \label{eq:mean_field_approximated_energy}
            &=\sum_{i=1}^{N_\text{spin}}h_ix_i+\sum_{i=1}^{N_\text{spin}}\sum_{j<i}Q_{i,j}x_ix_j
        \end{align}
        \quad\par hold. From \textbf{Eq.~\ref{eq:KL_difinition}$\sim$Eq.~\ref{eq:mean_field_approximated_energy}}, we obtain \textbf{Eq.~\ref{eq:KL_divergence}}. \par
        \quad Thus, the proposition has been proven.
\end{proof}

\subsection{Annealed Mean Field Descent (AMFD)}
\begin{figure}[b]
    \centering
    \includegraphics[width=0.92\linewidth]{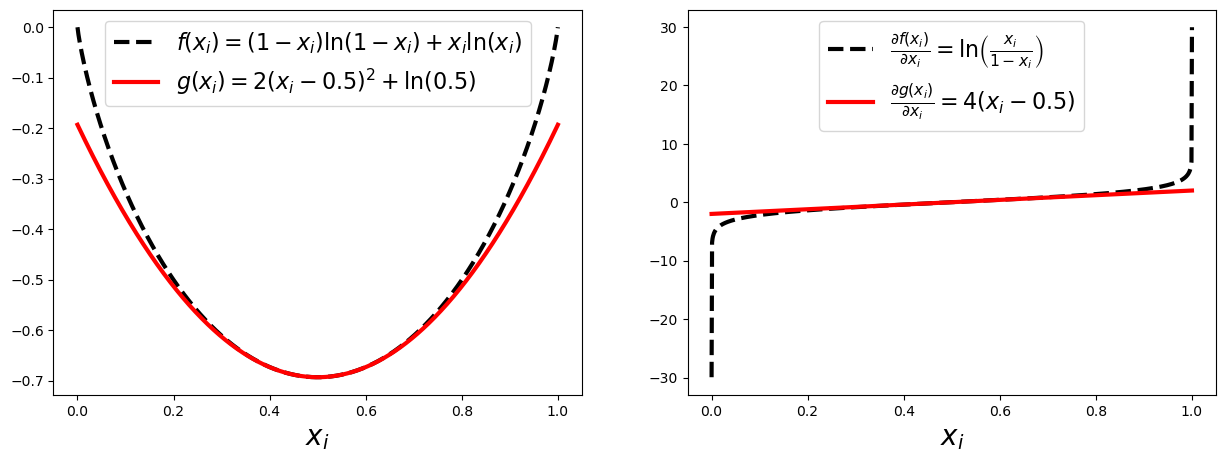}
    \caption{mf-entropy and its derivative.
Left/right: mf-entropy and derivative; black/red: original and Taylor approximation.}
    \label{fig:approximated_entropy}
\end{figure}
The theoretical analysis of the MFA reveals that the approximate distribution expressed in \textbf{Eq.~\ref{eq:approximate_distribution}} is a candidate probability distribution such that the KL divergence is minimized.
However, it is only the candidate of the minimum solution, not necessarily the minimum solution. Therefore, it is conceivable that the accuracy can be further improved by performing global minimization.
Then, we consider a minimization method using the gradient of the KL divergence to $\boldsymbol{x}$. The gradient can be calculated as follows:

\begin{align}
    \label{eq:KL_partial_by_x}
    &\frac{\partial}{\partial x_i}D_\text{KL}\left(P_\text{MF}\left(\boldsymbol{s}\right) \mid P_\text{C}\left(\boldsymbol{s}\right)\right) \notag\\
    &=\frac{\partial}{\partial x_i}\sum_{i=1}^{N_\text{spin}}\sum_{s_i=0}^{1}p_i(s_i)\ln\left(p_i(s_i)\right)
    + \frac{1}{T}\left(h_i + \sum_{j=1}^{N_\text{spin}}Q_{i,j}x_j\right) \notag\\
    &=\ln\left(\frac{x_i}{1-x_i}\right) +\frac{1}{T}\left(h_i + \sum_{j=1}^{N_\text{spin}}Q_{i,j}x_j\right), \\
    & \because \sum_{i=1}^{N_\text{spin}}\sum_{s_i=0}^{1}p_i(s_i)\ln\left(p_i(s_i)\right) \notag\\
    &= \sum_{i=1}^{N_\text{spin}} \left[p_i(0)\ln\left(p_i(0)\right) + p_i(1)\ln\left(p_i(1)\right)\right] \notag\\
    &= \sum_{i=1}^{N_\text{spin}}\left[\left(1-p_i(1)\right)\ln\left(1-p_i(1)\right) + p_i(1)\ln\left(p_i(1)\right)\right] \notag \\
    & = \sum_{i=1}^{N_\text{spin}}\left[\left(1-x_i\right)\ln\left(1-x_i\right) + x_i\ln\left(x_i\right)\right] .
\end{align}
The first term in \textbf{Eq.~\ref{eq:KL_partial_by_x}} implies the gradient of mf-entropy, and the second term implies the mean field (in other words, the gradient of mf-energy).

Here, since the second term in \textbf{Eq.~\ref{eq:KL_partial_by_x}} diverges as the temperature $T$ approaches $0$, our algorithm multiply \textbf{Eq.~\ref{eq:KL_partial_by_x}} by $T$ to prevent overflow. In addition, the gradient of mf-entropy (i.e., first term in \textbf{Eq.~\ref{eq:KL_partial_by_x}}) diverges as $x_i$ approaches $0$ or $1$, so our algorithm approximate mf-entropy to the second order using \textit{Taylor} series expansion at $x_i=0.5$:
\begin{align}\label{eq:Tayler_expansion}
    &\sum_{i=1}^{N_\text{spin}}\left[\left(1-x_i\right)\ln\left(1-x_i\right) + x_i\ln\left(x_i\right)\right]
    \notag \\\simeq &\sum_{i=1}^{N_\text{spin}}\left[2(x_i-0.5)^2+\ln(0.5)\right],
\end{align}
which is visualized in \textbf{Fig.~\ref{fig:approximated_entropy}}---although a more accurate approximation with a higher order is possible, this level of approximation exhibits sufficient performance. 
Also, when $x_i$ becomes $0$ or $1$ during simulation, we ignore the mean field term (i.e., the second term of \textbf{Eq.~\ref{eq:KL_partial_by_x}}) because the gradient of the mf-entropy term (i.e., the first term of \textbf{Eq.~\ref{eq:KL_partial_by_x}}) becomes a significantly larger ratio than the mean field term as $x_i$ approaches $0$ or $1$. 
Thus, our algorithm updates the mf-spin variables using the KL divergence multiplied by $T$:
\begin{align}\label{eq:KL_divergence_multiplied_by_T}
    &T\frac{\partial}{\partial x_i}D_\text{KL}\left(P_\text{MF}\left(\boldsymbol{s}\right) \mid P_\text{C}\left(\boldsymbol{s}\right)\right) \notag \\
    % &= T\ln\left(\frac{x_i}{1-x_i}\right) + h_i + \sum_{j=1}^{N_\text{spin}}Q_{i,j}x_j\\
    &\simeq\begin{cases}
    4T(x_i - 0.5) +  h_i + \sum_{j=1}^{N_\text{spin}}Q_{i,j}x_j, &(0 < x_i <1)\\
    4T(x_i - 0.5), &(x_i=0 \lor x_i=1)
    \end{cases}.
\end{align}
As a side note, $4T$ in \textbf{Eq.~\ref{eq:KL_divergence_multiplied_by_T}} can be lumped together as a new parameter $T'$, so the coefficient of $4$ is ignored in our algorithm.
In this way, scaling the gradient of the KL divergence with temperature not only prevents overflow at low temperatures but also has the secondary effect of increasing the step size at higher temperatures, making it easier to escape local minima.

Furthermore, unlike simple gradient descent, our algorithm incorporates an acceleration strategy that is expected to enhance convergence and facilitate escape from local solutions like Nesterov's Accelerated Gradient (NAG) method \cite{Nesterov1983AMF}. 
The slight differences from NAG are (1) the damping of the inertia term is not taken into account because the damping occurs when applying the constraints, and (2) the gradient calculation at the forward point is only applied to mf-energy and not to mf-entropy because the variables fluctuate violently and induce oscillations between 0 and 1. 

Also, in all experiments of this paper, the interaction matrix $\boldsymbol{Q}$ and external magnetic field vector $\boldsymbol{h}$ are normalized to make it easy to set parameters: 
\begin{equation}\label{eq:normalized_QUBO}
    \begin{cases}
        \boldsymbol{Q}_\text{norm} = \frac{\boldsymbol{Q}}{\sqrt{\frac{1}{N_\text{spin}}\sum_{i=1}^{N\text{spin}}\left(h_i^2 + \sum_{j=1}^{N_\text{spin}}Q_{i,j}^2\right)}}\\
        \boldsymbol{h}_\text{norm} = \frac{\boldsymbol{h}}{\sqrt{\frac{1}{N_\text{spin}}\sum_{i=1}^{N\text{spin}}\left(h_i^2 + \sum_{j=1}^{N_\text{spin}}Q_{i,j}^2\right)}}
    \end{cases}.
\end{equation}

\begin{figure*}[t]
    \centering
    \includegraphics[width=0.9\linewidth]{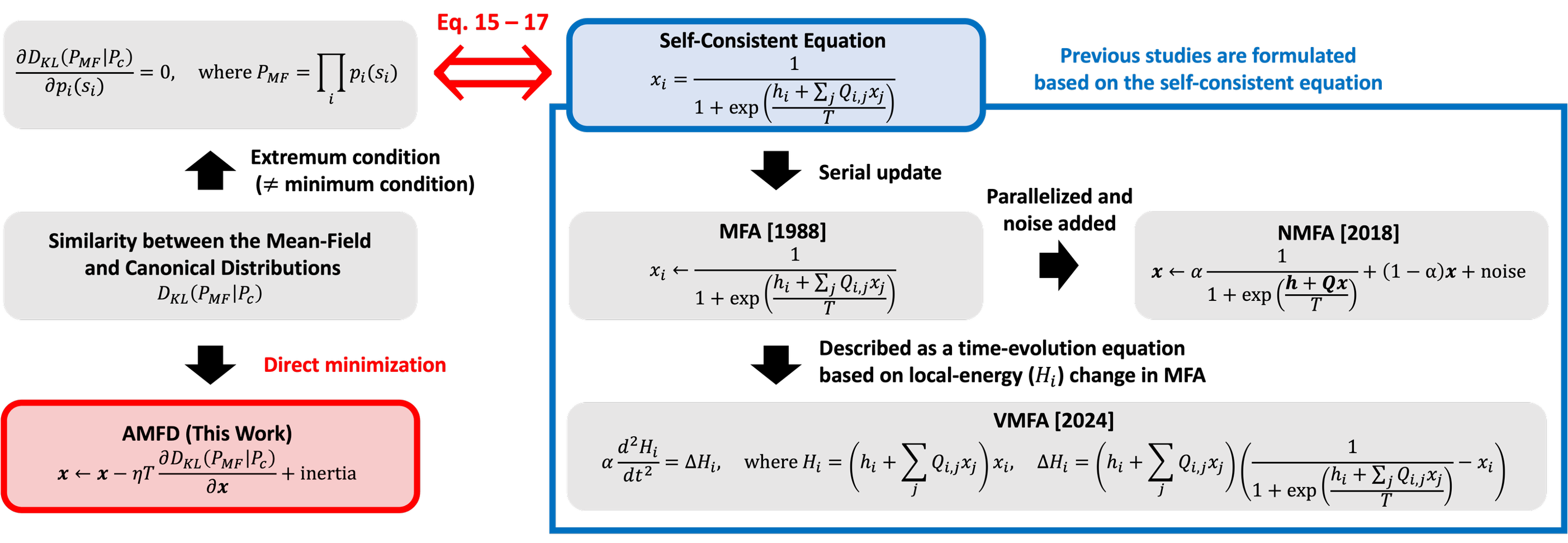}
    \caption{Positioning of this work. Previous studies update variables based on the extremum condition of the KL divergence between the canonical and mean-field distributions, whereas this work directly minimizes the KL divergence itself.}
    \label{fig:overview_of_AMFD}
\end{figure*}

Based on the above, we propose Annealed Mean Field Descent (AMFD) as shown in \textbf{Alg.~\ref{Alg:AMFD}}, and the overview of this work is displayed in \textbf{Fig.~\ref{fig:overview_of_AMFD}}.

\begin{algorithm}[h]
\caption{Annealed Mean Field Descent (AMFD)} \label{Alg:AMFD}
\normalsize
\begin{algorithmic}[1]
\Require $\bm{x}(-1)$, $\bm{x}(0)$, $\bm{Q}$, $\bm{h}$, $N_\text{step}$, $T(t)$, $\eta$, $\zeta$
\Ensure  $\bm{x}(N_\text{step})$
\State Normalization
\State $\bm{Q_\text{norm}} \gets \frac{\bm{Q}}{\sqrt{\frac{1}{N_\text{spin}}\sum_{i=1}^{N\text{spin}}\left(h_i^2 + \sum_{j=1}^{N_\text{spin}}Q_{i,j}^2\right)}}$
\State $\bm{h_\text{norm}} \gets \frac{\bm{h}}{\sqrt{\frac{1}{N_\text{spin}}\sum_{i=1}^{N\text{spin}}\left(h_i^2 + \sum_{j=1}^{N_\text{spin}}Q_{i,j}^2\right)}}$

\ForTo{$t = 1$}{$N_\text{step}$}
    
    \State Forward Point
    \State $\boldsymbol{x}_\text{fwd}(t) \gets \boldsymbol{x}(t-1) + \zeta\left(\boldsymbol{x}(t-1) - \boldsymbol{x}(t-2) \right)$
    \State Gradient of the MF-Entropy
    \State $\boldsymbol{F}(t) \gets T(t)\left(\boldsymbol{x}(t-1)-\boldsymbol{0.5}\right)$
    \State Mean Field at the Forward Point
    \State $\boldsymbol{{\Phi}}(t) \gets \boldsymbol{h_\text{norm}} + \boldsymbol{Q_\text{norm}}\boldsymbol{x}_\text{fwd}(t)$
    \ForToPara{$i = 1$}{$N_\text{spin}$}
        \State Descent with Acceleration
        \State ${x_i}(t) \gets 2{x_i}(t-1) - {x_i}(t-2) - \eta {F_i}(t)$ 
        \If {$0 < x_i(t-1) < 1$} 
    \State ${x_i}(t) \gets {x_i}(t)  - \eta {\Phi_i}(t)$ 
    \EndIf
        \State Apply Constraints
        \If {$x_i(t) \geq 1$} $x_i(t) \gets 1$ 
        \EndIf
        \If{$x_i(t) \leq 0$} $x_i(t) \gets 0$ 
        \EndIf 

    \EndForToPara
\EndForTo
\State Transition to High Probability State (Post-processing)
\ForToPara{$i = 1$}{$N_\text{spin}$}
    \If {$x_i(N_\text{step}) \geq 0.5$}
         \State  $x_i(N_\text{step}) \gets 1$  
    \Else 
        \State $x_i(N_\text{step}) \gets 0$ 
    \EndIf 
\EndForToPara
\end{algorithmic}
\end{algorithm}

Although the proposed method does not necessarily reach the global optimum solution, experimental results show that it reaches the optimum or near-optimum solution regardless of the type or instance of the COP. 

\section{Evaluation}\label{sec:evaluation}

This section evaluates AMFD compared to SOTA parallel SAs for QUBO (NMFA\cite{king18_NMFA}, bSB/dSB\cite{Goto_21_bSB}, RPA\cite{Kawamura23_RPA}, and VMFA\cite{Kuroki_VMFA}), learning-based methods (PI-GNN \cite{Schuetz2022_PIGNN}, RUN-CSP \cite{Toenshoff2019RUN_CSP}, QRF-GNN \cite{pugacheva2024QRF_GNN}, and VCM \cite{Chen2024VCM}), and Gurobi\cite{gurobi}, a SOTA mathematical optimization solver that integrates various elaborate methods and the branch-and-cut method.
\red{Our implementation is available at \href{https://github.com/kyo-kuroki/AMFD}{https://github.com/kyo-kuroki/AMFD}.}

\subsection{Experimental Settings}\label{subsec:experimental_setting}
We experimented with five typical COPs—MCP, MISP, TSP, QAP, and GCP—as well as general QUBO instances.
For each problem, we used standard benchmark datasets: Gset~\cite{Ye03_Gset} for MCP, the DIMACS benchmark set~\cite{Mascia_MaximumClique_Webpage} for MISP, TSPLIB~\cite{TSPlib, Reinelt_TSPLIB95} for TSP, QAPLIB~\cite{Burkard1997_QAPLIB_Webpage} for QAP, Graph Coloring Instances~\cite{Trick_COLOR_instances, DIMACS_GCP} for GCP, and general QUBO instances from ORLIB~\cite{orlib}.

 Each COP can be formulated into MILP, MIQP, and QUBO problems following the approach described in~\cite{glover2019tutorial, Lucas14_Formulation, DeSimone1990_MCP_as_MILP, Zhang2013_QAP_as_MILP, Gavish1978TheTS_TSP_as_MILP, MIS_MILP_sym15111979, Coll2002_GCP_MILP, Muramatsu2003ANS_MCP_as_MIQP, MIS_penalty, GCP_penalty}, and the detailed formulations are provided in supplementary materials.
Since Gurobi can handle MILP, MIQP, and QUBO formats, experiments of Gurobi are conducted for all formats to ensure maximum performance, rather than relying solely on the QUBO formulation. 
As the evaluation indicator, the gap is calculated using the Best Known Solution (BKS) and the obtained solution (SOL) by
\begin{equation}\label{eq:accuracy}
    \text{gap [\%]}= 100\cdot\left|\frac{\text{BKS} - \text{SOL}}{\text{BKS}}\right|.
\end{equation}

First, to compare the performance with existing parallel-SA-based QUBO solvers, we experimented the transition of the evaluation function (i.e., energy function) with respect to $N_\text{step}$ for each COP and each algorithm. For each COP, approximately 2000 to 2500 variable instances were selected. For all problems, the average and best values of the evaluation function $H$ on 1000 runs with $N_\text{step}$ swept were plotted for each algorithm. For all solvers, we used the best parameters selected through grid search. (See supplementary materials for the details.)
Also, it should be noted that all the parallel-SA-based QUBO solvers have almost the same computational cost per step, as the per-step complexity is dominated by the matrix–vector multiplication. 

Next, following the experimental setups in \cite{Schuetz2022_PIGNN, pugacheva2024QRF_GNN, Chen2024VCM}, we compared our approach with existing learning-based methods using seven MCP instances and ten general QUBO instances (bqp2500.1--bqp2500.10 from \cite{orlib}).
The results of the learning-based baselines were taken directly from the original studies: those for the MCP instances from \cite{Schuetz2022_PIGNN, pugacheva2024QRF_GNN}, and those for the general QUBO instances from \cite{Chen2024VCM}. The BKS values for the bqp2500 instances were obtained from \cite{WANG20105758_bqp_bestknown}.

Finally, we compared the performance of our method with that of Gurobi, conducting experiments on 40 instances for each COP.
 In the experiments, the AMFD was implemented and run on a GPU: NVIDIA GeForce RTX 4090 to exploit its high parallelism, and Gurobi was executed on two CPUs: AMD EPYC 7702 64-Core Processor using 256 threads---the version of Gurobi is 11.0.3.
Since Gurobi does not require hyperparameter tuning, we compare it against AMFD including its parameter tuning time, in order to ensure a fair comparison.

\begin{figure}[t]
    \centering
    \includegraphics[width=1\linewidth]{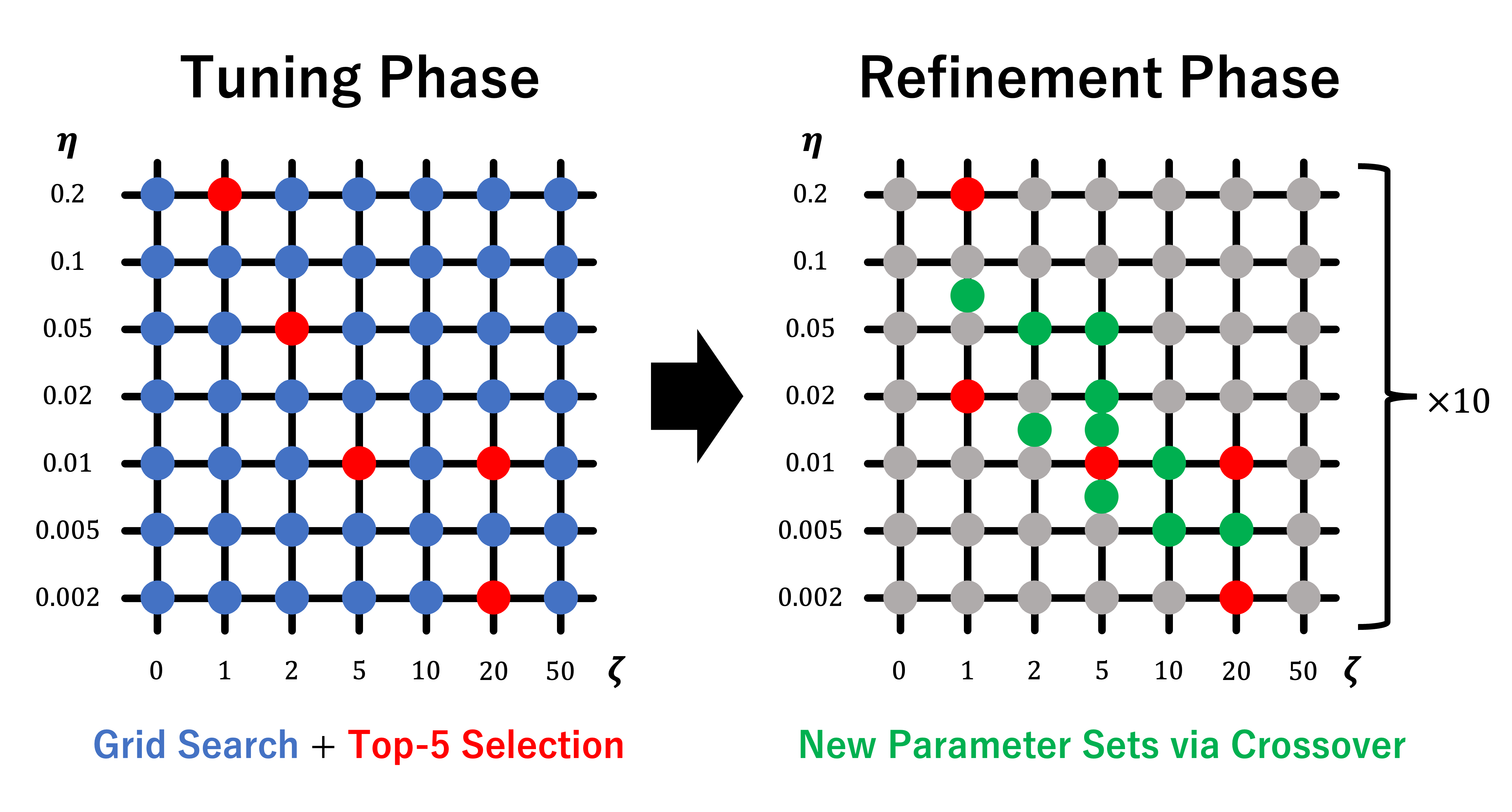}
    \caption{Parameter tuning using grid search and a genetic algorithm-based strategy.}
    \label{fig:parameter_tuning}
\end{figure}

For parameter tuning, as illustrated in \textbf{Fig.~\ref{fig:parameter_tuning}}, we adopt a hybrid strategy combining grid search and genetic algorithm. Specifically, we construct a grid of 49 combinations: $\eta \in \{0.002, 0.005, 0.01, 0.02, 0.05, 0.1, 0.2\}$, $\zeta \in \{0, 1, 2, 5, 10, 20, 50\}$, with $T_\text{init} = 0.35$ and $T_\text{end} = 0.001$. Each of the 49 combinations is run with a single replica, resulting in 49 parallel grid search executions.
These are executed in parallel on a GPU with $N_\text{step} = 2 N_\text{spin}$.
This tuning phase is performed only once in principle. However, when the formulation and the solution obtained through tuning indicate that certain variables can be eliminated, the tuning phase is repeated until no further unnecessary variables remain in the original problem.
 For example, in the case of the GCP, it means that the maximum number of available colors does not need to exceed that in the solution obtained during the tuning phase.

After the tuning phase, we extract the top-5 solutions (parent selection) and form new parameter sets by computing centroids of all pairwise $(\eta, \zeta)$ combinations (crossover): 
\[
\left\{ \left( \tfrac{1}{2} \eta_i + \tfrac{1}{2} \eta_j,\ \tfrac{1}{2} \zeta_i + \tfrac{1}{2} \zeta_j \right) \,\middle|\, 1 \le i < j \le 5 \right\},
\]
resulting in 10 additional parameter sets.
Together with the original top-5, we obtain 15 parameter sets for the post-tuning phase.
In this phase, we perform 10 replicas for each of the 15 parameter sets, resulting in 150 parallel executions. For each parameter set, we run experiments using five different values of $N_\text{step}$: $2 \times, 5 \times, 10 \times, 20 \times$, and $50 \times N_\text{spin}$. For each setting, we evaluate both the quality of the obtained solutions and the total computation time, including the tuning phase.

For Gurobi, the time limit was set to twice the maximum execution time of AMFD (i.e., the sum of tuning time and runtime with $N_\text{step} = 50 N_\text{spin}$), providing Gurobi with ample time for solution search. Whenever a solution was found, we recorded its value and the time at which it was obtained.

In experiments, the temperature schedule $T(t)$ for the AMFD was set according to 
\begin{equation}\label{eq:temperature_schedule}
    T(t)=T_\text{init} - \frac{T_\text{init}- T_\text{fin}}{N_\text{step}-1}(t-1),
\end{equation}
where $T_\text{init}$ represents the initial temperature, and $T_\text{fin}$ represents the final temperature. 
Also, the initial values of the mf-spin $\boldsymbol{x}(-1)$ were randomly initialized within the range $[0, 1]$, and $\boldsymbol{x}(0)$ was set as $\boldsymbol{x}(0) = \boldsymbol{x}(-1) - \eta(\boldsymbol{x}(-1) - 0.5)$. 
Other hyperparameters were determined by grid search.  
Details on the parameter settings are provided in supplementary materials.

\subsection{Experimental Results}\label{subsec:experimental_results}
\subsubsection{Comparison with SOTA parallel SAs for QUBO}\label{subsubsec:comparison_with_QUBO_solvers}
\begin{figure*}[t]
    \centering
    \includegraphics[width=0.7\linewidth]{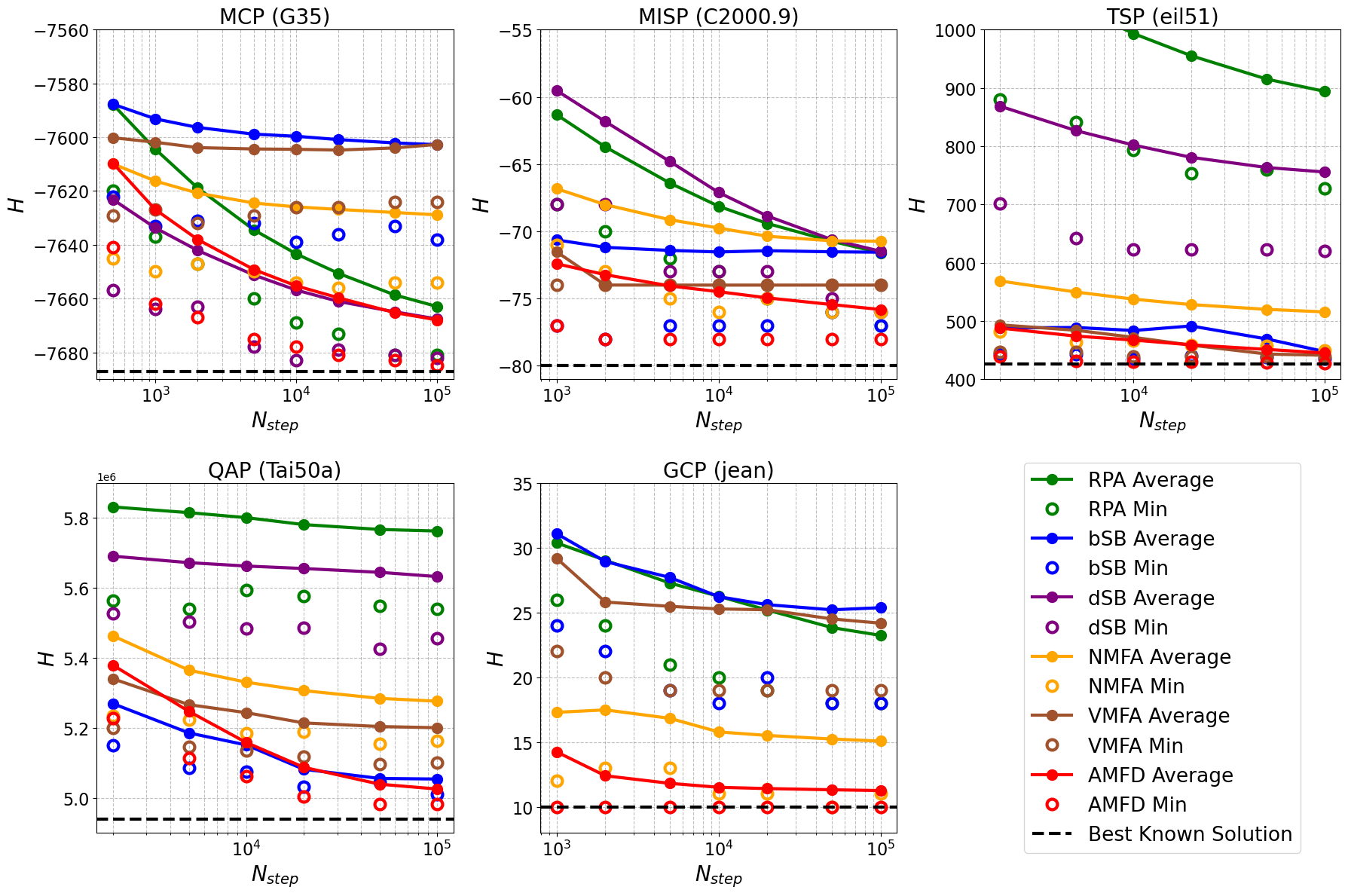}
    \caption{Comparison of parallel-SA-based QUBO solvers on a MCP, a MISP, a TSP, a QAP, and a GCP.}
    \label{fig:All_comparison}
\end{figure*}
\textbf{Fig.~\ref{fig:All_comparison}} exhibits the experimental results compared with SOTA parallel SAs for QUBO on five types of COPs.
AMFD consistently achieves high-quality solutions for all problems. On the other hand, the solution quality obtained by existing QUBO solvers is highly dependent on the type of COP---namely, they fail to achieve high-quality solutions for some types of COPs. 
For example, RPA and dSB are struggling to obtain high-quality solutions in QAP and TSP, while bSB and VMFA are struggling in MCP and GCP.
Therefore, it can be said that AMFD is more adaptable to a variety of COPs, which is a very valuable property for a general-purpose solver.

\vspace{1\baselineskip}
\subsubsection{Comparison with learning-based methods for QUBO}\label{subsubsec:comparison_with_learning_based_method}

\begin{table}[t]
\caption{Results of AMFD and learning-based baselines on four MCP benchmarks.}
\label{tab:comparison_with_learning_based_methods}
\centering
\footnotesize
\begin{tabular}{crrrrr}
\hline
Instance & \multicolumn{1}{c}{BKS} & \multicolumn{1}{c}{PI-GNN} & \multicolumn{1}{c}{RUN-CSP} & \multicolumn{1}{c}{QRF-GNN} & \multicolumn{1}{c}{AMFD} \\ \hline\hline
G14      & 3064                    & 3026                       & 2943                        & 3058                        & \textbf{3063}            \\
G15      & 3050                    & 2990                       & 2928                        & \textbf{3049}               & 3048                     \\
G22      & 13359                   & 13181                      & 13028                       & 13340                       & \textbf{13359}           \\
G49      & 6000                    & 5918                       & \textbf{6000}               & \textbf{6000}               & \textbf{6000}            \\
G50      & 5880                    & 5820                       & \textbf{5880}               & \textbf{5880}               & \textbf{5880}            \\
G55      & 10294                   & 10138                      & 10116                       & 10282                       & \textbf{10293}           \\
G70      & 9591                    & 9421                       & 9319                        & \textbf{9559}               & 9538                     \\ \hline
\end{tabular}
\end{table}

\begin{table}[t]
\caption{Average gap [\%] of AMFD and learning-based methods on ten general QUBO instances (bqp2500).}
\label{tab:evaluation_on_bqp2500}
\centering
\small
\begin{tabular}{crrr}
\hline
                  & \multicolumn{1}{c}{PI-GNN} & \multicolumn{1}{c}{VCM-BGF-HB} & \multicolumn{1}{c}{AMFD} \\ \hline
Avg. Gap {[}\%{]} & 1.689                      & 0.012                          & \textbf{0.000}           \\
Avg. Time {[}s{]} & 44                         & 19                             & \textbf{2.0}               \\ \hline
\end{tabular}
\end{table}

\textbf{Tab.~\ref{tab:comparison_with_learning_based_methods}} and \textbf{Tab.~\ref{tab:evaluation_on_bqp2500}} present the comparison results with the learning-based methods.
Specifically, \textbf{Tab.~\ref{tab:comparison_with_learning_based_methods}} shows the best solutions obtained by each method on seven MCP instances, while \textbf{Tab.~\ref{tab:evaluation_on_bqp2500}} reports the average gaps across ten general QUBO instances.
As shown in the tables, AMFD consistently provides higher-quality solutions than the learning-based methods across most problem instances, demonstrating its strong optimization capability.
Especially, for all bqp2500 instances, AMFD successfully reaches the BKS, which highlights the effectiveness of our approach.
\red{Furthermore, it is noteworthy that the proposed method operates within a substantially shorter runtime than these learning-based methods, and also avoids the need for additional network design or COP-specific training.}

\vspace{1\baselineskip}
\subsubsection{Comparison with Gurobi}\label{subsubsec:comparison_with_gurobi}

\begin{table}[t]
\caption{Average gap [\%] of AMFD and Gurobi over 40 instances across five types of COPs.}
\label{tab:average_gap}
\centering
\footnotesize
\begin{tabular}{crrrrr}
\hline
               & \multicolumn{1}{c}{MCP} & \multicolumn{1}{c}{MISP} & \multicolumn{1}{c}{TSP} & \multicolumn{1}{c}{QAP} & \multicolumn{1}{c}{GCP} \\ \hline\hline
Gurobi (MILP)  & 11.95                   & 11.55                    & \textbf{0.00}           & 11.42                   & 93.44                   \\
Gurobi (MIQP)  & 15.48                   & 19.29                    & 414.95                  & 11.48                   & 21.36                   \\
Gurobi (QUBO)  & 1.53                    & 2.94                     & 72.06                   & 10.97                   & 26.04                   \\
AMFD (QUBO)    & \textbf{0.17}           & \textbf{0.96}            & 1.56                    & \textbf{1.09}           & \textbf{2.47}           \\
to-AMFD (QUBO) & 1.05                    & 5.68                     & 10.47                   & 6.28                    & 8.47                    \\ \hline
\end{tabular}
\end{table}

\begin{figure}[h]
    \centering
    \includegraphics[width=0.9\linewidth]{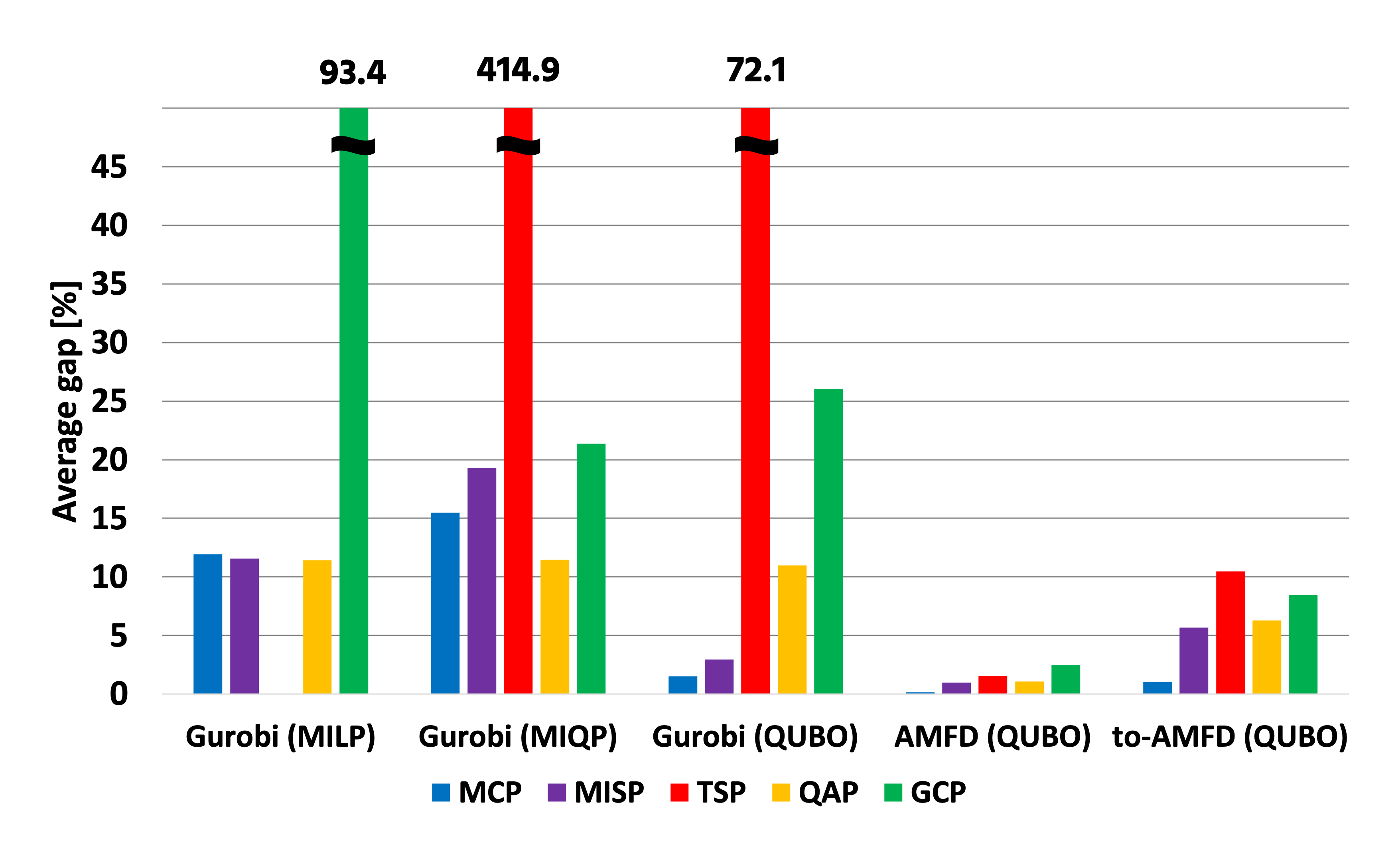}
    \caption{Average gap [\%] comparison of AMFD and Gurobi.}
    \label{fig:average_gap}
\end{figure}

\begin{table*}[h]
\caption{Evaluation of Gurobi and AMFD on MCPs.}
\label{tab:MCP_performance}
\centering
\footnotesize
\begin{tabular}{crrrrrrrrrrr}
\hline
\multirow{2}{*}{Instance} & \multicolumn{1}{c}{\multirow{2}{*}{BKS}} & \multicolumn{2}{c}{Gurobi (MILP)}                          & \multicolumn{2}{c}{Gurobi (MIQP)}                          & \multicolumn{2}{c}{Gurobi (QUBO)}                          & \multicolumn{2}{c}{AMFD (QUBO)}                            & \multicolumn{2}{c}{to-AMFD (QUBO)}                         \\
                          & \multicolumn{1}{c}{}                     & \multicolumn{1}{c}{SOL} & \multicolumn{1}{c}{Time {[}s{]}} & \multicolumn{1}{c}{SOL} & \multicolumn{1}{c}{Time {[}s{]}} & \multicolumn{1}{c}{SOL} & \multicolumn{1}{c}{Time {[}s{]}} & \multicolumn{1}{c}{SOL} & \multicolumn{1}{c}{Time {[}s{]}} & \multicolumn{1}{c}{SOL} & \multicolumn{1}{c}{Time {[}s{]}} \\ \hline\hline
G1                        & 11624                                    & 10306                   & 3.63                             & N/A                     & \textless{}15.84                 & 11598                   & 14.09                            & \textbf{11624}          & 2.44                             & 11623                   & 0.85                             \\
G2                        & 11620                                    & 10234                   & 3.67                             & N/A                     & \textless{}15.93                 & 11599                   & 9.12                             & \textbf{11620}          & 2.26                             & 11617                   & 0.92                             \\
G3                        & 11622                                    & 10271                   & 9.95                             & N/A                     & \textless{}16.36                 & 11612                   & 7.68                             & \textbf{11622}          & 0.99                             & \textbf{11622}          & 0.99                             \\
\hline
\end{tabular}
\end{table*}

% Please add the following required packages to your document preamble:
% \usepackage{multirow}
\begin{table*}[t]
\caption{Evaluation of Gurobi and AMFD on MISPs.}
\label{tab:MISP_performance}
\centering
\footnotesize
\begin{tabular}{crrrrrrrrrrr}
\hline
\multirow{2}{*}{Instance} & \multicolumn{1}{c}{\multirow{2}{*}{BKS}} & \multicolumn{2}{c}{Gurobi (MILP)}                          & \multicolumn{2}{c}{Gurobi (MIQP)}                          & \multicolumn{2}{c}{Gurobi (QUBO)}                          & \multicolumn{2}{c}{AMFD (QUBO)}                            & \multicolumn{2}{c}{to-AMFD (QUBO)}                         \\
                          & \multicolumn{1}{c}{}                     & \multicolumn{1}{c}{SOL} & \multicolumn{1}{c}{Time {[}s{]}} & \multicolumn{1}{c}{SOL} & \multicolumn{1}{c}{Time {[}s{]}} & \multicolumn{1}{c}{SOL} & \multicolumn{1}{c}{Time {[}s{]}} & \multicolumn{1}{c}{SOL} & \multicolumn{1}{c}{Time {[}s{]}} & \multicolumn{1}{c}{SOL} & \multicolumn{1}{c}{Time {[}s{]}} \\ \hline\hline
brock200\_1               & 21                                       & \textbf{21}             & 4.67                             & \textbf{21}             & 4.37                             & \textbf{21}             & 1.56                             & \textbf{21}             & 3.06                             & 20                      & 1.48                             \\
brock400\_1               & 27                                       & 21                      & 2.59                             & 22                      & 7.11                             & \textbf{25}             & 7.63                             & \textbf{25}             & 3.62                             & 24                      & 1.55                             \\
brock800\_1               & 23                                       & 17                      & 18.6                             & 17                      & 16.72                            & 20                      & 8.19                             & \textbf{21}             & 3.96                             & 20                      & 1.69                             \\
\hline
\end{tabular}
\end{table*}

% Please add the following required packages to your document preamble:
% \usepackage{multirow}
\begin{table*}[t]
\caption{Evaluation of Gurobi and AMFD on TSPs.}
\label{tab:TSP_performance}
\centering
\footnotesize
\begin{tabular}{crrrrrrrrrrr}
\hline
\multirow{2}{*}{Instance} & \multicolumn{1}{c}{\multirow{2}{*}{BKS}} & \multicolumn{2}{c}{Gurobi (MILP)}                          & \multicolumn{2}{c}{Gurobi (MIQP)}                          & \multicolumn{2}{c}{Gurobi (QUBO)}                          & \multicolumn{2}{c}{AMFD (QUBO)}                            & \multicolumn{2}{c}{to-AMFD (QUBO)}                         \\
                          & \multicolumn{1}{c}{}                     & \multicolumn{1}{c}{SOL} & \multicolumn{1}{c}{Time {[}s{]}} & \multicolumn{1}{c}{SOL} & \multicolumn{1}{c}{Time {[}s{]}} & \multicolumn{1}{c}{SOL} & \multicolumn{1}{c}{Time {[}s{]}} & \multicolumn{1}{c}{SOL} & \multicolumn{1}{c}{Time {[}s{]}} & \multicolumn{1}{c}{SOL} & \multicolumn{1}{c}{Time {[}s{]}} \\ \hline\hline
att48                     & 10628                                    & \textbf{10628}          & 4.28                             & 15474                   & 40.82                            & 16730                   & 9.91                             & 10729                   & 9.31                             & 11278                   & 1.54                             \\
bays29                    & 2020                                     & \textbf{2020}           & 0.18                             & 2354                    & 16.41                            & 2164                    & 16.46                            & 2026                    & 7.99                             & 2034                    & 1.01                             \\
berlin52                  & 7542                                     & \textbf{7542}           & 0.41                             & 11239                   & 40.42                            & 10273                   & 48.23                            & 7553                    & 23.36                            & 8160                    & 1.72                             \\
\hline
\end{tabular}
\end{table*}

% Please add the following required packages to your document preamble:
% \usepackage{multirow}
\begin{table*}[t]
\caption{Evaluation of Gurobi and AMFD on QAPs.}
\label{tab:QAP_performance}
\centering
\footnotesize
\begin{tabular}{crrrrrrrrrrr}
\hline
\multirow{2}{*}{Instance} & \multicolumn{1}{c}{\multirow{2}{*}{BKS}} & \multicolumn{2}{c}{Gurobi (MILP)}                          & \multicolumn{2}{c}{Gurobi (MIQP)}                          & \multicolumn{2}{c}{Gurobi (QUBO)}                          & \multicolumn{2}{c}{AMFD (QUBO)}                            & \multicolumn{2}{c}{to-AMFD (QUBO)}                         \\
                          & \multicolumn{1}{c}{}                     & \multicolumn{1}{c}{SOL} & \multicolumn{1}{c}{Time {[}s{]}} & \multicolumn{1}{c}{SOL} & \multicolumn{1}{c}{Time {[}s{]}} & \multicolumn{1}{c}{SOL} & \multicolumn{1}{c}{Time {[}s{]}} & \multicolumn{1}{c}{SOL} & \multicolumn{1}{c}{Time {[}s{]}} & \multicolumn{1}{c}{SOL} & \multicolumn{1}{c}{Time {[}s{]}} \\ \hline\hline
bur26a                    & 5426670                                  & 5611826                 & 14.37                            & 5528569                 & 12.12                            & 5500764                 & 9.96                             & \textbf{5484383}        & 3.86                             & 5644287                 & 0.99                             \\
chr25a                    & 3796                                     & 5000                    & 13.5                             & 4492                    & 11.22                            & 5130                    & 5.91                             & \textbf{3832}           & 6.62                             & 6256                    & 0.97                             \\
esc32a                    & 130                                      & 150                     & 20.3                             & 152                     & 18.03                            & 162                     & 19.71                            & \textbf{130}            & 9.92                             & 166                     & 1.15                             \\
\hline
\end{tabular}
\end{table*}

% Please add the following required packages to your document preamble:
% \usepackage{multirow}
\begin{table*}[t]
\caption{Evaluation of Gurobi and AMFD on GCPs.}
\label{tab:GCP_performance}
\centering
\footnotesize
\begin{tabular}{crrrrrrrrrrr}
\hline
\multirow{2}{*}{Instance} & \multicolumn{1}{c}{\multirow{2}{*}{BKS}} & \multicolumn{2}{c}{Gurobi (MILP)}                          & \multicolumn{2}{c}{Gurobi (MIQP)}                          & \multicolumn{2}{c}{Gurobi (QUBO)}                          & \multicolumn{2}{c}{AMFD (QUBO)}                            & \multicolumn{2}{c}{to-AMFD (QUBO)}                         \\
                          & \multicolumn{1}{c}{}                     & \multicolumn{1}{c}{SOL} & \multicolumn{1}{c}{Time {[}s{]}} & \multicolumn{1}{c}{SOL} & \multicolumn{1}{c}{Time {[}s{]}} & \multicolumn{1}{c}{SOL} & \multicolumn{1}{c}{Time {[}s{]}} & \multicolumn{1}{c}{SOL} & \multicolumn{1}{c}{Time {[}s{]}} & \multicolumn{1}{c}{SOL} & \multicolumn{1}{c}{Time {[}s{]}} \\ \hline\hline
anna                      & 11                                       & \textbf{11}             & 3.23                             & \textbf{11}             & 1.16                             & 12                      & 1.84                             & \textbf{11}             & 8.15                             & \textbf{11}             & 8.15                             \\
david                     & 11                                       & \textbf{11}             & 1.65                             & \textbf{11}             & 1.23                             & 12                      & 1.51                             & \textbf{11}             & 8.25                             & \textbf{11}             & 8.25                             \\
DSJC125.1                 & 5                                        & 6                       & 17.51                            & 6                       & 5.96                             & 7                       & 15.63                            & \textbf{5}              & 11.16                            & 6                       & 7.74                             \\
\hline
\end{tabular}
\end{table*}

\textbf{Fig.~\ref{fig:average_gap}} and \textbf{Tab.~\ref{tab:average_gap}} show the average gaps for each problem type---MCP, MISP, TSP, QAP, and GCP---computed over 40 instances each, based on the best solutions obtained. Note that instances for which no feasible solution was found are excluded from the average calculations.
In addition, \textbf{Tab.~\ref{tab:MCP_performance},~\ref{tab:MISP_performance},~\ref{tab:TSP_performance},~\ref{tab:QAP_performance},~\ref{tab:GCP_performance}} present the experimental results comparing AMFD and Gurobi across a total of 200 instances.
In each instance row of \textbf{Tab.~\ref{tab:MCP_performance}--\ref{tab:GCP_performance}}, we report, for each method, the best solution obtained and the time taken to reach it. 
N/A or values in parentheses indicate cases where no feasible solution satisfying all constraints was obtained; the value inside the parentheses represents the QUBO objective value of the infeasible solution.
For AMFD, the reported time corresponds to the total duration of both the tuning and refinement phases. On the other hand, to-AMFD (tuning-only AMFD) denotes the setting where only the tuning phase is performed; for this setting, we report the best solution obtained during the tuning phase together with the corresponding tuning time. Due to space limitation, we show only the best solutions and their first achievement times; more detailed results are in supplementary materials.

Across all problem types evaluated, AMFD consistently delivered high-quality solutions with small optimality gaps, often outperforming Gurobi regardless of the problem formulation (\textbf{Tab.~\ref{tab:average_gap}}, \textbf{Fig.~\ref{fig:average_gap}}). Specifically, it achieved an average gap below $3\%$ for all COPs considered—MCP, MISP, TSP, QAP, and GCP. While Gurobi excelled in certain cases—particularly with well-suited formulations or specific graph structures—its performance tended to vary more by instance, whereas AMFD demonstrated stable effectiveness across diverse settings.

A closer look at each COP reveals the following:
\begin{itemize}
    \item \textbf{MCP:} AMFD achieved $\leq 1\%$ gaps for all instances. Even with only the tuning phase (referred to as to-AMFD), it also reached near-optimality quickly. Gurobi excelled on toroidal graphs (e.g., \texttt{G32}, \texttt{G34}), indicating a structural advantage.
    \item \textbf{MISP:} AMFD again led overall with an average gap within $1\%$. Gurobi’s QUBO formulation performed comparably, but MILP/MIQP versions struggled on large instances ($\geq 1000$ nodes), highlighting the scalability of the QUBO form.
    \item \textbf{TSP:} Gurobi’s MILP solved all instances to optimality, but its QUBO results were poor (average gap $>$72\%). AMFD remained competitive (1.56\% average gap), reinforcing its strength in QUBO solving.
    \item \textbf{QAP:} AMFD outperformed all others except on \texttt{tai64c} and \texttt{esc128}, where gaps were still small (0.27\% and 3.1\%). The average gap was 1.09\%, showing near-optimal performance.
    \item \textbf{GCP:} AMFD achieved an average gap of 2.47\% with consistent quality, while Gurobi’s performance varied greatly and failed on certain instances (e.g., \texttt{school1}), indicating higher robustness for AMFD.
\end{itemize}

\section{Discussion}

\subsection{Factors Contributing to AMFD’s High Performance.}
First, we discuss why AMFD performs better on many COPs than other SOTA parallel SAs. 
As shown in \textbf{Fig.~\ref{fig:overview_of_AMFD}}, AMFD directly minimizes the KL divergence to achieve a mean-field approximate distribution as accurate as possible to the true distribution, whereas the self-consistent equation (\textbf{Eq.~\ref{eq:self_consistent_equation}}) on which NMFA and VMFA are based exhibits the extremum condition for the KL divergence. Therefore, we believe that the solution quality of AMFD is often better than NMFA or VMFA because it is easier to realize a more accurate approximate distribution.
In addition, the reason why RPA, a MCMC-based parallel SA, does not outperform AMFD is thought to be that an enormous number of sampling iterations are required to approach the true distribution in the MCMC process, and as a result, it does not converge to a better distribution than the mean-field approximate distribution within a given time. In particular, RPA seems not to obtain good solutions in TSP, QAP, and GCP. This is likely because these COPs involve many functions that represent constraints, and when these constraints are violated, the energy rapidly increases, resulting in a steep and complex landscape. In such cases, once MCMC-based SA traps a local minimum of deep valleys, it tends to repeatedly sample around there, making it difficult to escape.
Also, the reason why bSB and dSB, quantum annealing inspired algorithms, cannot obtain high-quality solutions for some COPs is thought to be that they define the Hamiltonian using continuous variables, which does not accurately map the QUBO function. In contrast, although the mf-spin variables in AMFD themselves are continuous values, they are calculated based on the canonical distribution defined by the QUBO function composed of binary variables. Therefore, if the error due to the continuous value definition is larger than the error due to the mean-field approximation, AMFD is expected to give better results than bSB and dSB.

\subsection{AMFD as a General-Purpose Solver.}
Next, we discuss AMFD's performance as a general-purpose solver through comparisons with Gurobi. 
The purpose of the comparisons is not to provide a way of selecting the optimal solver based on which is superior between AMFD and Gurobi for each problem. Rather, the comparisons aim to observe which of the two solvers, as a general-purpose solver, can consistently deliver high-quality solutions regardless of the problem type or instance.
From this perspective, the overall experimental results indicate that AMFD consistently achieves higher-quality solutions than Gurobi across various problem types and instances, demonstrating its superior functionality as a general-purpose solver.
Surely, just as Gurobi with a MILP formulation is more efficient in solving TSP, some COPs may be more efficiently solved by Gurobi with a MILP or MIQP formulation---it can be inferred that such problems have a characteristic that makes it easy to narrow the lower and upper bounds when formulated as MILP or MIQP. However, it is generally difficult to introduce formulations that are easy to narrow the bounds, and also, the fact that solver performance can vary greatly depending on the formulation is not suitable for the democratization of combinatorial optimization.

\subsection{Further Potential.}\label{subsec:further_potential}
In contrast to Gurobi, which is a composite solver based on the branch-and-cut framework and incorporates a variety of heuristics, AMFD is built upon a single metaheuristic. Nevertheless, it is noteworthy that AMFD achieved superior average optimality gaps compared to Gurobi—whose time limit was set to twice the maximum execution time of AMFD—on many COPs.
 This observation suggests that integrating AMFD into a branch-and-cut framework, for example as a component for tightening upper bounds, could lead to further performance improvements in composite solvers.
Moreover, in the field of QUBO optimization, recent studies have explored decomposition-based approaches that partition a problem into smaller subproblems and subsequently integrate their solutions \cite{ISCA22_divide, Guerreschi2021SolvingQU_divide_and_conquer}. Combining AMFD with such decomposition techniques could further accelerate solution times and improve solution quality for large-scale instances.

Additionally, although AMFD is based on a mean-field approximation, the underlying framework could, in principle, be extended to more general variational approximations. For instance, the current approach assumes independence in the probability distribution of individual spins; this independence assumption could be generalized to joint distributions over spin pairs or other richer variational families, potentially improving accuracy. However, such extensions would inevitably increase the number of variational parameters, potentially reducing computational efficiency. Moreover, while the mean-field approximation allows explicit gradients of the KL divergence with respect to variational expectations, extended models do not retain this property.
 This necessitates additional efforts such as incorporating penalty functions for constraints or assuming specific functional forms for the variational distribution, making the approach more challenging. 
Nevertheless, there remains significant potential for improved accuracy, making the exploration and evaluation of these richer approximations a worthwhile direction for future research.

Furthermore, in this paper, the penalty coefficients for constraints in solving constrained COPs were determined statically based on previous studies and theoretical analysis. However, as in \cite{Huang2021QROSSQR, Roch_optimize_penalty_coeff}, dynamically optimizing these coefficients could potentially lead to further performance improvements.

\section{Conclusion}
This paper proposes AMFD, a method designed to efficiently solve COPs formulated in QUBO (a format capable of expressing various types of COPs), aiming to advance the democratization of combinatorial optimization. 
AMFD is designed to more accurately approximate the true canonical distribution, motivated by the observation that the mean-field approximation error of MFA and its variants can be further reduced.
Through extensive experiments, it has been demonstrated that AMFD provides a high-quality solution that is less dependent on problem types and instances compared to existing methods. 

Furthermore, there is potential for further performance improvements through branch-and-bound techniques employed by Gurobi, divide-and-conquer approaches, enriched variational distributions, and streamlined parameter optimization.
We believe that AMFD opens a new avenue toward the democratization of combinatorial optimization.

\clearpage
\bibliographystyle{IEEEtran}
\bibliography{reference} 
\begin{IEEEbiography}[{\includegraphics[width=1in,height=1.25in,clip,keepaspectratio]{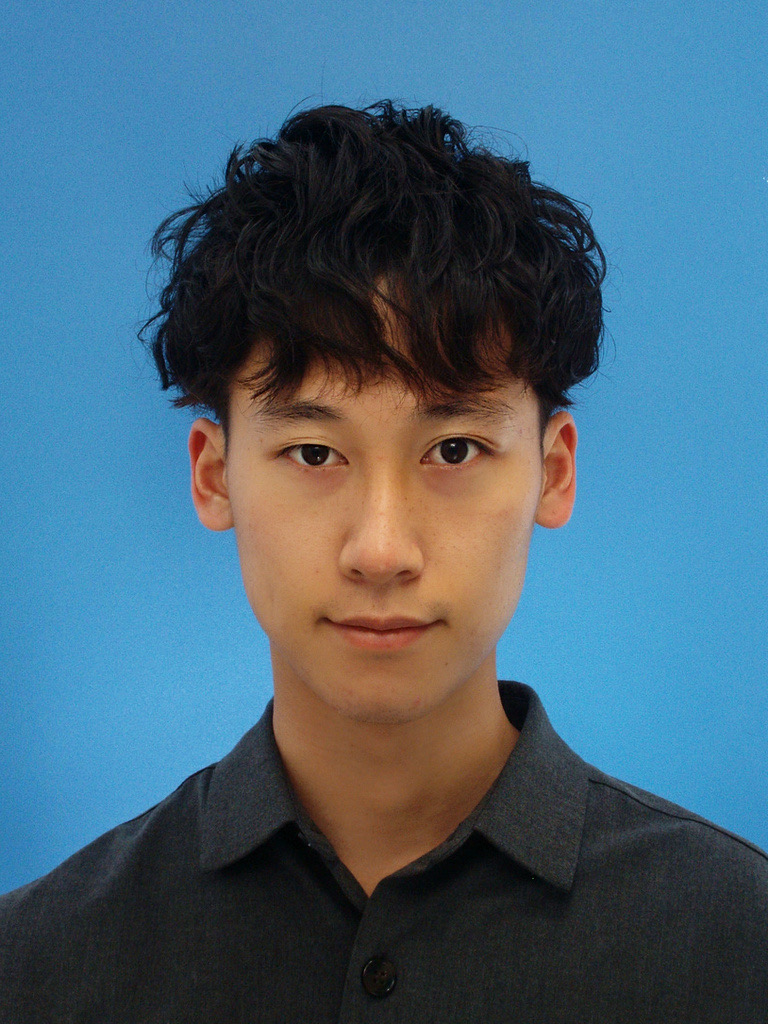}}]{Kyo Kuroki}
received the B.E. degree in engineering from Tohoku University and the M.E. degree in information and communications engineering from Institute of Science Tokyo (formerly Tokyo Institute of Technology). He is currently pursuing the Ph.D. degree in information and communications engineering at Institute of Science Tokyo. 

His research interests include evolutionary computation, combinatorial optimization, and efficient machine learning.
\end{IEEEbiography}

\begin{IEEEbiography}[{\includegraphics[width=1in,height=1.25in,clip,keepaspectratio]{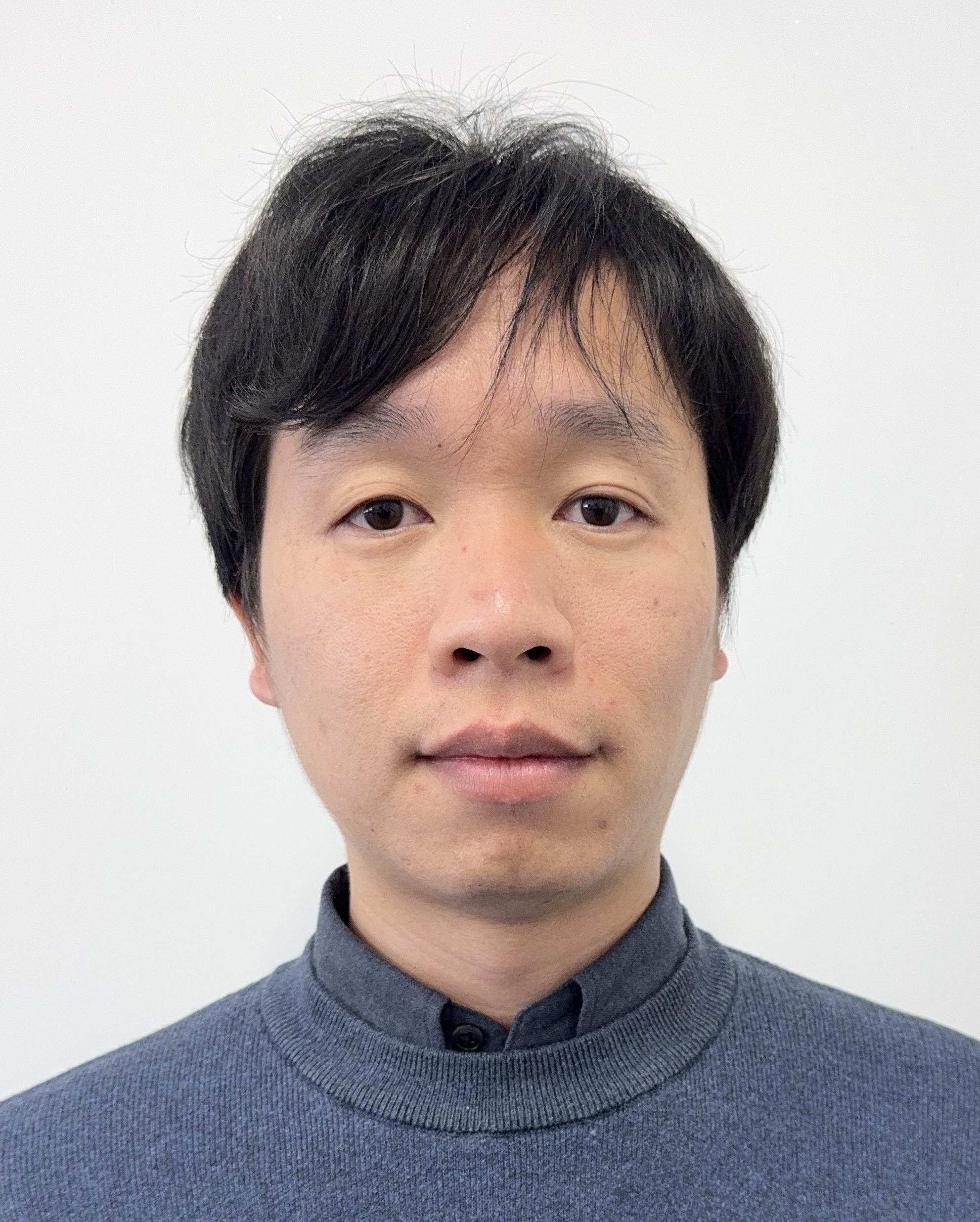}}]{Thiem Van Chu} (Member, IEEE) is an Associate Professor in the Department of Information and Communications Engineering at Institute of Science Tokyo. He received his Ph.D. in Engineering from Tokyo Institute of Technology in 2018. His research interests include computer architecture, systems, and machine learning, with recent emphasis on efficient algorithms and hardware for machine learning and sparse data processing. He has received several honors, including Best Paper Awards at the International Conference on Field-Programmable Technology (FPT) in 2021 and the International Symposium on Field-Programmable Custom Computing Machines (FCCM) in 2017.
\end{IEEEbiography}

\begin{IEEEbiography}[{\includegraphics[width=1in,height=1.25in,clip,keepaspectratio]{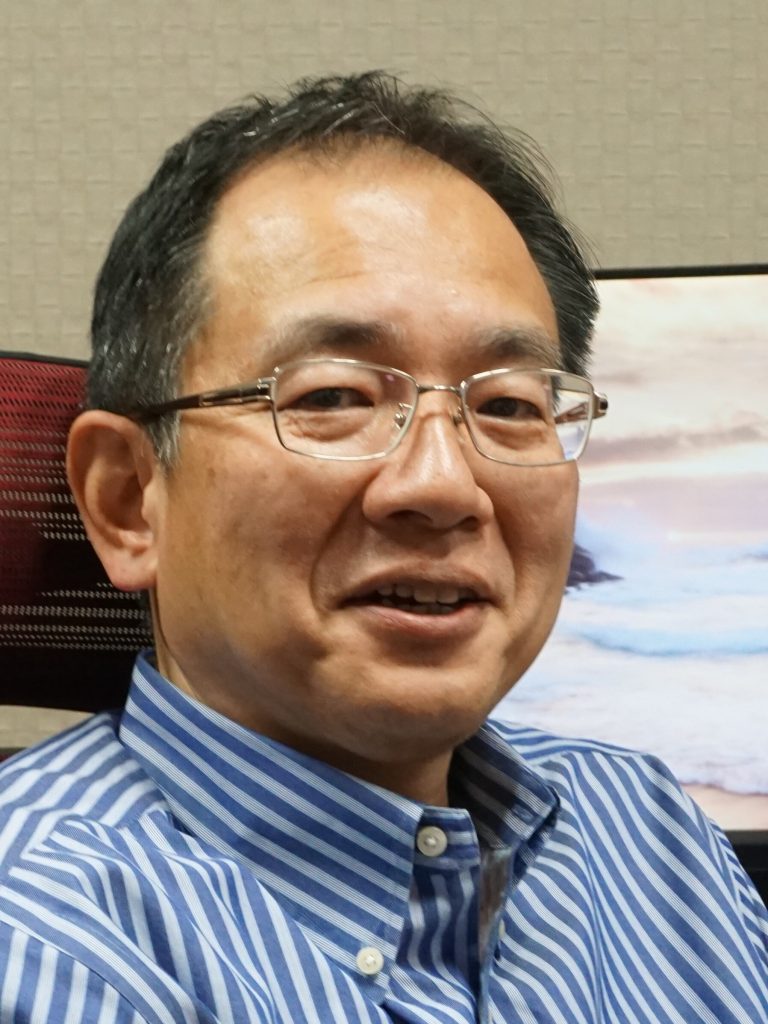}}]{Masato Motomura} (Fellow, IEEE) received BS, MS, and Dr. of Engineering from Kyoto University, Kyoto, Japan, in 1985, 1987, and 1996, respectively.
In 1987, he joined NEC Central Research Laboratories, Kawasaki, Japan. He was a visiting researcher at MIT Laboratory for Computer Science, Cambridge, USA, from 1991 to 1992. He led research and business development of the dynamically reconfigurable processor (DRP) he invented from 2001 to 2008.

In 2011,  he became a professor at Hokkaido University, Sapporo, Japan.
Later he became a professor at Tokyo Institute of Technology (TokyoTech), Yokohama, Japan, in 2019, where he established and has been leading the artificially intelligent computing (ArtIC) research unit.
Since 2011, he has been actively working on reconfigurable and parallel architectures for deep neural networks, machine learning, annealing machines, and general intelligent/domain-specific computing.
He received the IEEE JSSC Annual Best Paper Award in 1992, the IPSJ Annual Best Paper Award in 1999, and the IEICE Achievement Award in 2011, respectively.
He was also awarded  Ichimura Academic Award in 2022,
Yamasaki Award in 2023, and  MEXT Award in 2024, for his leadership in developing and productizing DRP technology (a series of DRP-based microcontroller products are now produced by Renesas Electronics), as well as the accumulation of AI-chip achievements in recent years.
He is a member of IEICE, IPSJ, JSAI, and EAJ.
He is a 2022 IEEE Fellow (SSCS) for contributions to memory-logic integration of reconfigurable chip architecture.
\end{IEEEbiography}

\begin{IEEEbiography}
[{\includegraphics[width=1in,height=1.25in,clip,keepaspectratio]{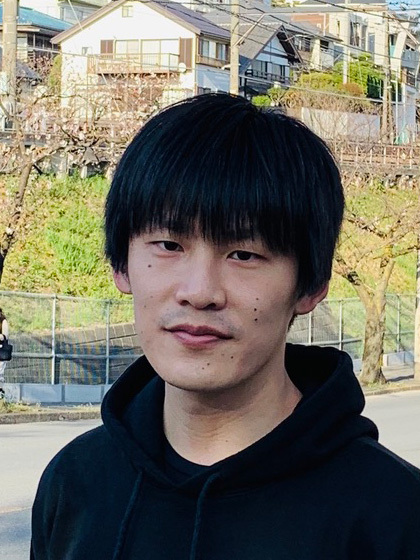}}]
{Kazushi Kawamura} (Member, IEEE) received the B.Eng., M.Eng., and Dr.Eng. degrees from Waseda University, in 2012, 2013, and 2016, respectively, all in computer science. From 2018 to 2019, he was an Assistant Professor with the Department of Communications and Computer Engineering, Waseda University. In 2020, he joined the Institute of Innovative Research, Tokyo Institute of Technology, as a Specially Appointed Assistant Professor. He is currently an Associate Professor with the Department of Electronic and Physical Systems, Waseda University. His research interests include intelligent processing systems, particularly those based on combinatorial optimization and machine learning. He is a member of ACM, IEICE, and IPSJ.
\end{IEEEbiography}

\clearpage
\pagestyle{empty} 
% \appendix
\section*{Supplementary Materials}
\renewcommand{\thefigure}{S.\arabic{figure}}
\renewcommand{\thetable}{S.\arabic{table}}
\renewcommand{\theequation}{S.\arabic{equation}}
\renewcommand{\thealgorithm}{S.\arabic{algorithm}}  % algorithm2e等を使っている場合

% カウンタのリセット
\setcounter{figure}{0}
\setcounter{table}{0}
\setcounter{equation}{0}
\setcounter{algorithm}{0}

\renewcommand{\thesection}{S.\arabic{section}}
\renewcommand{\thesubsection}{S.\arabic{subsection}}
\setcounter{section}{0}  % subsectionカウンタをリセット
\setcounter{subsection}{0}  % subsectionカウンタをリセット

\section{Formulation of COPs}\label{sec:cop_formulation}
\vspace{1\baselineskip}
\subsubsection{Maximum Cut Problem (MCP)}\label{subsec:MCP}
MCP is the problem that maximizes the weight of edges between different groups when dividing a given set of vertices into two subsets. We used Gset \cite{Ye03_Gset} as the datasets for the experiments.

Formulating MCP as QUBO, we obtain the following function:
\begin{equation}
\label{eq:MCP_QUBO}
    H = -\frac{1}{2}\sum_{i = 1}^{N_{\text{vertex}}}\sum_{j < i}E_{i,j}\left[1-(2s_i-1)(2s_j-1)\right],
\end{equation}
where $E_{i,j} \in \mathbb{R}$ denotes the weight of the edge connecting vertex $i$ and vertex $j$, $N_\text{vertex}$ represents the number of vertices, and $s_i \in \{0, 1\}$ is a variable that equals 1 if the vertex belongs to one group and 0 if it belongs to the other group after the split. The value of $H$ represents the total number of cuts. The negative sign is introduced to transform the problem into a minimization objective.

Also, formulating the MCP as a MILP \cite{DeSimone1990_MCP_as_MILP} leads to
\begin{align}\label{eq:MCP_MILP}
    \text{minimize}\quad &-\sum_{i=1}^{N_\text{vertex}}\sum_{j<i}E_{i,j}y_{i,j},\\
    \text{s.t. \quad} & s_i + s_j \ge y_{i,j}, \quad \forall (i,j) \in E, \\
    & 2 - s_i -s_j \ge y_{i,j}, \quad \forall (i,j) \in E, \\
    & s_i - s_j \le y_{i,j}, \quad \forall (i,j) \in E, \\
    & s_j - s_i \le y_{i,j}, \quad \forall (i,j) \in E, 
\end{align}
where $E = \{(i,j)|E_{i,j}\ne 0 \land j < i\}$ represents the set of edges in the graph and $y_{i,j}\in \{0, 1\}$ represents a variable that is $1$ when vertices $i$ and $j$ are included on different groups and $0$ when they are included on the same group.

In MCP formulation as MIQP, it is known that it can be formulated as a quadratic cone optimization problem\cite{Muramatsu2003ANS_MCP_as_MIQP} as shown below:
\begin{align}\label{eq:MCP_MIQP}
    \text{minimize} \quad &-\sum_{i=1}^{N_\text{vertex}}\sum_{j<i}E_{i,j}z_{i,j},\\
    \text{s.t. \quad} & (s_i + s_j - 1)^2 \le u_{i,j}, \quad \forall (i,j) \in E, \\
    & (s_i - s_j)^2 \le z_{i,j}, \quad \forall (i,j) \in E,\\
    & z_{i,j} + u_{i,j} = 1, \quad \forall (i,j) \in E, 
\end{align}
where $u_{i,j} \in \{0, 1\}$ represents the variable that is 1 when vertices $i$ and $j$ are included on the same group, and $z_{i,j}  \in \{0, 1\}$ represents the variable that is 1 when vertices $i$ and $j$ are included on different group.

\vspace{1\baselineskip}
\subsubsection{Maximum Independent Set Problem (MISP)}\label{subsec:MISP}
Given a graph consisting of vertices and edges, MISP is defined as the problem of finding the largest vertex set where each vertex is not connected by an edge. 
Since MISP is equivalent to the maximum clique problem when considering the complement graph of a given graph, our experiments used the complement graphs of the maximum clique problem datasets \cite{Mascia_MaximumClique_Webpage}.

The QUBO formulation for MISP is expressed by
\begin{equation}
    \label{eq:MISP_QUBO}
            H = -\sum_{i=1}^{N_{\text{vertex}}}s_{i} + A\sum_{(i,j)\in E}s_{i}s_{j},
\end{equation}
where $E$ represents the set of edges in the graph, $N_\text{vertex} \in \mathbb{N}$ represents the number of vertices, $s_i \in \{0,1\}$ represents $1$ when vertex $i$ is included in the independent set and $0$ when it is not, and $A \in \mathbb{R}$ is a constant coefficient. The first term---with the minus sign to attribute MISPs to minimization problems---refers to the size of the set, and the second term represents the penalty imposed for violating the constraint where no vertex in the set is connected to any other vertex by an edge.
% In our experiments, we set $A = 2$  to make it more dominant than the objective term.
{In our experiments, since the maximum change in the objective function value is $1$, we set the penalty coefficient $A = 2$ based on \cite{MIS_penalty} to ensure that the constraint term is more dominant.
}

The MILP formulation for MISP\cite{MIS_MILP_sym15111979} is expressed by
\begin{align}\label{eq:MISP_MILP}
    \text{minimize} \quad &-\sum_{i=1}^{N_\text{vertex}}s_i,\\
    \label{eq:MISP_MILP_constraint}
    \text{s.t. \quad} & s_i + s_j \le 1, \quad \forall (i,j) \in E,
\end{align}
where \textbf{Eq.~\ref{eq:MISP_MILP_constraint}} represents the constraint that vertices connected by edges are not included in the independent set.

Also, the MIQP formulation for MISP can be expressed by
\begin{align}\label{eq:MISP_MIQP}
    \text{minimize} \quad &-\sum_{i=1}^{N_\text{vertex}}s_i,\\
    \label{eq:MISP_MIQP_constraint}
    \text{s.t. \quad} & s_is_j = 0, \quad \forall (i,j) \in E.
\end{align}
As with the MILP formulation, \textbf{Eq.~\ref{eq:MISP_MIQP_constraint}} represents the constraint that vertices connected by edges are not included in the independent set.

\vspace{1\baselineskip}
\subsubsection{Traveling Salesman Problem (TSP)}\label{subsec:TSP}
Given a set of cities and a travel cost (e.g. distance) between each of the two cities, TSP is the problem of finding the tour that travels through all cities exactly once and returns to the starting city with the lowest total travel cost. 

Given a travel cost $D_{i,j}$ between city $i$ and city $j$, the variable $s_{i,k} \in \{0,1\}$ visiting city $i$ in $k$-th position, the QUBO formulation of the TSP is described as 
\begin{equation}\label{eq:TSP_QUBO}
    \begin{split}
    H =& \sum_{i=1}^{N_{\text{city}}-1}\sum_{j=1}^{N_{\text{city}}-1}\sum_{k=1}^{N_{\text{city}}-2}D_{i,j}s_{i,k}s_{j,(k+1)} \\
    &+ \sum_{i=1}^{N_{city}-1}D_{N_\text{city}, i}s_{i,1}+ \sum_{i=1}^{N_{city}-1}D_{i,N_\text{city}}s_{i,(N_\text{city}-1)} 
    \\&+ A\sum_{k=1}^{N_{\text{city}-1}}\left(1-\sum_{i=1}^{N_{\text{city}}-1}s_{i,k}\right)^2 
    \\&+ B\sum_{i=1}^{N_{\text{city}}-1}\left(1-\sum_{k=1}^{N_{\text{city}}-1}s_{i,k}\right)^2,
    \end{split}
\end{equation}
where $N_\text{city}\in \mathbb{N}$ represents the number of cities, $A \in \mathbb{R}$ and $B \in \mathbb{R}$ are constant coefficients---set to $A=B=\max_i \left\{\frac{\sum_{j}D_{i,j}}{\alpha\cdot(N_\text{city}-1)}\right\}$ which represents maximum average distance from city $i$. 
This is based on the assumption that, in the optimal solution, the inter-city distance does not exceed the maximum, over all \(i\), of the average distance from city \(i\). (Note that \(\alpha\) represents the density of the inter-city distance matrix \(\boldsymbol{D}\).)

The first, second, and third terms represent the total travel cost of the tour, the fourth term represents the penalty imposed for violating the constraint where only one city is visited at a time, and the fifth term represents the penalty imposed for violating the constraint where each city $i$ is visited only once. Note that the first and last city visited is set to the $N_\text{city}$-th city.

Next, we describe the MILP formulation of the TSP. Although various formulations of MILP for TSP have been proposed, in our experiments, we use the formulation \cite{Gavish1978TheTS_TSP_as_MILP} known to be efficient, as follows:
\begin{align}\label{eq:TSP_MILP}
    \text{minimize} \quad&\sum_{i=1}^{N_\text{city}}\sum_{j=1}^{N_\text{city}}D_{i,j}y_{i,j},\\
    \text{s.t. \quad} & \sum_{j\ne i}y_{i,j} = 1, \quad \forall i, \label{eq:TSP_MILP_const1}\\
    & \sum_{j\ne i}y_{j,i} = 1, \quad \forall i, \label{eq:TSP_MILP_const2}\\
    & \sum_{j=2}^{N_\text{city}}f_{1,j} = N_\text{city} - 1, \label{eq:TSP_MILP_const3}\\
    & \sum_{j\ne i}f_{j,i} - \sum_{j\ne i}f_{i,j} = 1, \quad \text{for} \quad i\ge 2,\label{eq:TSP_MILP_const4}\\
    & f_{1,i} \le (N_\text{city} - 1)y_{1, i} \quad \text{for} \quad i \ge 2,  \label{eq:TSP_MILP_const5}\\
    & f_{i,j} \le (N_\text{city} - 2)y_{i,j} \quad \forall i \ne j, i \ne 1, j\ne1
    , \label{eq:TSP_MILP_const6}
\end{align}
where $f_{i,j}\in \{r| r\in \mathbb{R} \land r \ge 0\}$ is a variable representing the order of the tour, with $(N_\text{city}-1)$ quantities flowing from the starting point to the next point to be visited, with the quantity decreasing by one for each visit, $y_{i,j}\in \{0, 1\}$ is a variable that represents $1$ when passing from city $i$ to city $j$ and $0$ when not passing.
\textbf{Eq.~\ref{eq:TSP_MILP_const1}, \ref{eq:TSP_MILP_const2}} are constraints that prevent movement from one point to multiple points or from multiple points to one point, and \textbf{Eq.~\ref{eq:TSP_MILP_const3}, \ref{eq:TSP_MILP_const4}} imply constraints that $(N_\text{city}-1)$ quantities flow from the first point to the second point and that the difference between inflow and outflow is $1$ at all points except the first and last city visited.
\textbf{Eq.~\ref{eq:TSP_MILP_const5}, \ref{eq:TSP_MILP_const6}} denote the upper limits of flow.
Note that city $1$ is designated as the first and last city visited in this formulation.

Also, the MIQP formulation for TSP can be expressed by
\begin{align}\label{eq:TSP_MIQP}
    \text{minimize}\quad & \sum_{i=1}^{N_{\text{city}}-1}\sum_{j=1}^{N_{\text{city}}-1}\sum_{k=1}^{N_{\text{city}}-2}D_{i,j}s_{i,k}s_{j,(k+1)} \notag \\
    &+ \sum_{i=1}^{N_{city}-1}D_{N_\text{city}, i}s_{i,1}+ \sum_{i=1}^{N_{city}-1}D_{i,N_\text{city}}s_{i,(N_\text{city}-1)} ,
    \\ \text{s.t. \quad}&  \sum_{i=1}^{N_{\text{city}}-1}s_{i,k} = 1, \quad \forall k,\label{eq:TSP_MIQP_const1}
    \\&\sum_{k=1}^{N_{\text{city}}-1}s_{i,k} = 1, \quad \forall i,\label{eq:TSP_MIQP_const2}
\end{align}
where \textbf{Eq.~\ref{eq:TSP_MIQP_const1}, \ref{eq:TSP_MIQP_const2}} represent constraints that you cannot visit multiple cities at the same time and that you can only visit the same city once.

\vspace{1\baselineskip}
\subsubsection{Quadratic Assignment Problem (QAP)}\label{subsec:QAP}

Assuming that the flow amount between factory $i$ and factory $j$ is $F_{i,j}$ and the distance between city $i$ and city $j$ is $D_{i,j}$, the problem is to find the allocation that minimizes the total transportation cost---which is expressed as the product of the flow amount and distance.

Formulating the QAP as QUBO, QUBO function can be expressed as
\begin{equation}\label{eq:QAP_QUBO}
    \begin{split}
    H =& \sum_{i=1}^{N_{\text{city}}}\sum_{j=1}^{N_{\text{city}}}\sum_{k=1}^{N_{\text{city}}}\sum_{l=1}^{N_{\text{city}}}F_{i,k}D_{j,l}s_{i,j}s_{k,l} 
    \\&+ A\sum_{i=1}^{N_{\text{city}}}\left(1-\sum_{j=1}^{N_{\text{city}}}s_{i,j}\right)^2 
    + B\sum_{j=1}^{N_{\text{city}}}\left(1-\sum_{i=1}^{N_{\text{city}}}s_{i,j}\right)^2,
    \end{split}
\end{equation}
where $N_\text{city} \in \mathbb{N}$ represents the number of cities and the number of factories---the number of factories coincides with the number of cities---, $s_{i,j} \in \{0, 1\}$ represents a variable that becomes $1$ when factory $i$ is allocated to city $j$ and becomes $0$ when not allocated.
The first term implies as low a transportation cost as possible, the second term implies a constraint that each factory can be assigned to only one city, and the third term implies a constraint that only one factory can be assigned to each city.
We set constant coefficient $A$ and $B$ to $A=B=\max_{(i,j)} \left\{\sum_{k=1}^{N_\text{city}}\sum_{l=1}^{N_\text{city}}\frac{F_{i,k}D_{j,l}}{\alpha\cdot(N_\text{city}-1)}\right\}$, which represents maximum average transportation cost of allocating factory $i$ to city $j$.
{ Similar to the TSP case, this is based on the assumption that the transportation cost between optimally assigned factories is not greater than the maximum, over all \((i, j)\), of the average transportation cost when factory \(i\) is assigned to city \(j\). Here, \(\alpha\) denotes the combined density of the distance matrix \(\boldsymbol{D}\) and the flow matrix \(\boldsymbol{F}\), representing the ratio of nonzero elements in both matrices to the total number of elements in \(\boldsymbol{D}\) and \(\boldsymbol{F}\).
}

Next, considering $w_{i,j} \in \mathbb{R}$ as a variable representing the transportation cost of assigning factory $i$ to city $j$, the MILP formulation for QAP \cite{Zhang2013_QAP_as_MILP} can be expressed as
\begin{align}\label{eq:QAP_MILP}
    \text{minimize} \quad& \sum_{i=1}^{N_{\text{city}}}\sum_{j=1}^{N_{\text{city}}}w_{i,j} ,
    \\ \text{s.t. \quad} & \sum_{j=1}^{N_{\text{city}}}s_{i,j} = 1, \quad \forall i,\label{eq:QAP_MILP_const1}
    \\ &\sum_{i=1}^{N_{\text{city}}}s_{i,j} = 1, \quad \forall j, \label{eq:QAP_MILP_const2}
    \\ & \left(\sum_{k=1}^{N_{\text{city}}}\sum_{l=1}^{N_{\text{city}}}F_{i,k}D_{j,l}\right)(s_{i,j} - 1) \notag \\&+ \sum_{k=1}^{N_{\text{city}}}\sum_{l=1}^{N_{\text{city}}}F_{i,k}D_{j,l}s_{k,l} \le w_{i,j}, \quad \forall i, j,\label{eq:QAP_MILP_const3}
\end{align}
where \textbf{Eq.~\ref{eq:QAP_MILP_const1}, \ref{eq:QAP_MILP_const2}} represent the constraints that only one factory is allocated to one city and one factory is allocated to only one city. 
\textbf{Eq.~\ref{eq:QAP_MILP_const3}} represents the constraint that $w_{i,j}$ becomes the total transportation cost between factory $i$ and the other factories when factory $i$ is assigned to city $j$.

Also, the QAP can be formulated as MIQP, as shown below:
\begin{align}\label{eq:QAP_MIQP}
    \text{minimize} \quad& \sum_{i=1}^{N_{\text{city}}}\sum_{j=1}^{N_{\text{city}}}\sum_{k=1}^{N_{\text{city}}}\sum_{l=1}^{N_{\text{city}}}F_{i,k}D_{j,l}s_{i,j}s_{k,l} ,
    \\\text{s.t. \quad} & \sum_{j=1}^{N_{\text{city}}}s_{i,j} = 1, \quad \forall i,\label{eq:QAP_MIQP_const1}
    \\ &\sum_{i=1}^{N_{\text{city}}}s_{i,j} = 1, \quad \forall j,\label{eq:QAP_MIQP_const2}
\end{align}
where \textbf{Eq.~\ref{eq:QAP_MIQP_const1}, \ref{eq:QAP_MIQP_const2}} mean the constraints that only one factory is allocated to one city and one factory is allocated to only one city.

\vspace{1\baselineskip}
\subsubsection{Graph Coloring Problem (GCP)}\label{subsec:GCP}
GCP is the problem of painting vertices connected by edges with different colors from each other for a given graph and finding the minimum number of colors that can be painted in such a way.

Using $E$ representing the set of edges, the variable $s_{i,k} \in \{0,1\}$ representing $1$ when vertex $i$ is painted with color $k$ and $0$ when it is not, and the variable $y_{k}\in \{0,1\}$ that is $1$ when at least color $k$ is used and $0$ when not used, the QUBO formulation of GCP can be described as 
\begin{align}\label{eq:GCP_QUBO}
        H &= \sum_{k=1}^{N_{\text{color}}}y_k + A\sum_{k=1}^{N_\text{color}}\left[(1 - y_k)\sum_{i=1}^{N_\text{vertex}}s_{i,k}\right] \notag 
        \\ &+ B\sum_{i=1}^{N_{\text{vertex}}}\left(1-\sum_{k=1}^{N_{\text{color}}}s_{i,k}\right)^2  + C\sum_{k=1}^{N_{\text{color}}}\sum_{(i,j)\in E}s_{i,k}s_{j,k} , 
\end{align}
where $N_\text{color}$ represents the maximum number of colors that can be used---since it is proved that any graph can be colored with the maximum degree of the vertex plus one color, we set $N_\text{color} = \max\left\{\text{degree}_i\right\} + 1$ in our experiments---, $N_\text{vertex}$ represents the number of vertices, $A, B, C \in \mathbb{R}$ is a constant coefficient---{we set the penalty coefficients $A = B = C = 2$ in our experiments, based on \cite{GCP_penalty}, to ensure that the constraint terms are more dominant than the objective term, given that the maximum change in the objective function value is 1.}
The first term means the objective term and the second term represents the penalty for violating the constraint that $y_k$ must be 1 if at least color $k$ is used.
Note that when the color $k$ is not used, the second term is zero regardless of the value of $y_k$, but thanks to the first term, this formulation is sufficient because $H$ is lower when $y_k$ is $0$ than $y_k$ is $1$.
The third term represents the penalty for violating the constraint that only one color can be painted on the same vertex, and the fourth term represents the penalty for violating the constraint that vertices connected by an edge cannot be painted the same color.

Next, the MILP formulation of GCP\cite{Coll2002_GCP_MILP} can be described as
\begin{align}\label{eq:GCP_MILP}
        \text{minimize}\quad & \sum_{k=1}^{N_{\text{color}}}y_k, \\
        \text{s.t.\quad} &  \sum_{k=1}^{N_{\text{color}}}s_{i,k} = 1, \quad \forall i,
        \label{eq:GCP_MILP_const1}
        \\ &s_{i,k} + s_{j,k} \le y_k , \quad\forall k, (i,j)\in E, 
        \label{eq:GCP_MILP_const2}
\end{align}
where \textbf{Eq.~\ref{eq:GCP_MILP_const1}} represents the constraint that a vertex must not be painted with multiple colors and \textbf{Eq.~\ref{eq:GCP_MILP_const2}} represents the constraint that vertices connected with an edge must not be painted with the same color and that $y_k$ must be $1$ when color $k$ is used and must be $0$ when not used.

Also, the MIQP formulation of GCP can be described as
\begin{align}\label{eq:GCP_MIQP}
        \text{minimize}\quad & \sum_{k=1}^{N_{\text{color}}}y_k, \\
        \text{s.t.\quad} &  \sum_{k=1}^{N_{\text{color}}}s_{i,k} = 1, \quad \forall i,\label{eq:GCP_MIQP_const1}
        \\ &s_{i,k}s_{j,k} =  0 , \quad\forall k, (i,j)\in E, \label{eq:GCP_MIQP_const2}\\
        &(1 - y_k)\sum_{i=1}^{N_\text{vertex}}s_{i,k} = 0, \quad \forall k,\label{eq:GCP_MIQP_const3}\\
        &y_k \le \sum_{i=1}^{N_\text{vertex}}s_{i,k}, \quad \forall k\label{eq:GCP_MIQP_const4}
\end{align}
where \textbf{Eq.~\ref{eq:GCP_MIQP_const1}} represents the constraint that a vertex must not be painted with multiple colors, \textbf{Eq.~\ref{eq:GCP_MIQP_const2}} represents the constraint that vertices connected with an edge must not be painted with the same color, \textbf{Eq.~\ref{eq:GCP_MIQP_const3}} represents the constraint that $y_k$ must be $1$ when color $k$ is used, and \textbf{Eq.~\ref{eq:GCP_MIQP_const4}} represents the constraint that $y_k$ must be $0$ when color $k$ is not used.

\section{Parallel SAs’ Hyperparameters}

\begin{table}[h]
    \centering
    \caption{RPA Parameters for \textbf{Sec.\ref{subsubsec:comparison_with_QUBO_solvers}}}
    \label{tab:RPA_parameters}
    % \fontsize{10pt}{0.6cm}\selectfont
    \begin{tabular}{cccc}
\hline
{Instance} & {$\varepsilon$} & $T_\text{init}$ & {$T_\text{fin}$} \\ \hline\hline
MCP (G35)         & 0.5              & $0.5M$                    & $0.5M/N_\text{spin}$                    \\
MISP (C2000.9)    & 0.8              & $M$                    & $M/N_\text{spin}$                    \\
TSP (eil51)       & 0.2              & $M$                 & $0.5M/N_\text{spin}$                  \\
QAP (tai50a)      & 0.3              & $M$                   & $2M/N_\text{spin}$                    \\
GCP (jean)        & 0.1              & $0.5M$                 & $M/N_\text{spin}$                    \\ \hline
\end{tabular}
\end{table}

\begin{table}[h]
    \centering
    \caption{bSB Parameters for \textbf{Sec.\ref{subsubsec:comparison_with_QUBO_solvers}}}
    \label{tab:bSB_parameters}
    % \fontsize{10pt}{0.6cm}\selectfont
\begin{tabular}{cccc}
\hline
{Instance} & {$\Delta t$} & {$a_0$} & {$c_0$} \\ \hline\hline
MCP (G35)         & 1               & 1              & 0.5/G            \\
MISP (C2000.9)    & 0.25            & 0.5            & 1/G              \\
TSP (eil51)       & 0.25            & 0.5            & 1/G            \\
QAP (tai50a)      & 0.25            & 0.5            & 1/G              \\
GCP (jean)        & 0.75            & 0.5            & 0.5/G            \\ \hline
\end{tabular}
\end{table}

\begin{table}[h]
    \centering
    \caption{dSB Parameters for \textbf{Sec.\ref{subsubsec:comparison_with_QUBO_solvers}}}
    \label{tab:dSB_parameters}
    % \fontsize{10pt}{0.6cm}\selectfont
\begin{tabular}{cccc}
\hline
{Instance} & {$\Delta t$} & {$a_0$} & {$c_0$} \\ \hline\hline
MCP (G35)         & 0.75               & 1              & 1/G            \\
MISP (C2000.9)    & 0.25            & 0.5            & 0.5/G              \\
TSP (eil51)       & 0.25            & 0.5            & 1/G            \\
QAP (tai50a)      & 0.25            & 0.5            & 1/G              \\
GCP (jean)        & 0.25            & 0.5            & 0.5/G            \\ \hline
\end{tabular}
\end{table}

\begin{table}[h]
    \centering
    \caption{NMFA Parameters for \textbf{Sec.\ref{subsubsec:comparison_with_QUBO_solvers}}}
    \label{tab:NMFA_parameters}
    % \fontsize{10pt}{0.6cm}\selectfont
\begin{tabular}{ccccc}
\hline
Instance       & \multicolumn{1}{c}{$\alpha$} & \multicolumn{1}{c}{$\sigma$} & \multicolumn{1}{c}{$T_\text{init}$} & \multicolumn{1}{c}{$T_\text{fin}$} \\ \hline\hline
MCP (G35)      & 0.2                           & 0.1                         & 1                                    & 0.001                                   \\
MISP (C2000.9) & 0.2                           & 0.001                         & 0.01                                 & 0.001                                   \\
TSP (eil51)    & 0.3                           & 0.0001                        & 0.005                                 & 0                                   \\
QAP (tai50a)   & 0.2                           & 0.0001                        & 0.005                                 & 0.001                                   \\
GCP (jean)     & 0.1                           & 0.0001                        & 0.005                                & 0.0001                                   \\ \hline
\end{tabular}
\end{table}

\begin{table}[h]
    \centering
    \caption{VMFA Parameters for \textbf{Sec.\ref{subsubsec:comparison_with_QUBO_solvers}}}
    \label{tab:VMFA_parameters}
    % \fontsize{10pt}{0.6cm}\selectfont
\begin{tabular}{cccc}
\hline
Instance       & $\eta$ & $T_\text{init}$ & $T_\text{fin}$ \\ \hline\hline
MCP (G35)      & 0.5   & $V$               & $V/N_\text{spin}$            \\
MISP (C2000.9) & 0.25    & $V$             & $V/N_\text{spin}$            \\
TSP (eil51)    & 0.25    & $V$               & $V/N_\text{spin}$              \\
QAP (tai50a)   & 0.05   & $0.5V$               & $2V/N_\text{spin}$            \\
GCP (jean)     & 0.1    & $V$             & $0.5V/N_\text{spin}$            \\ \hline
\end{tabular}
\end{table}

\begin{table}[h]
    \centering
    \caption{AMFD Parameters for \textbf{Sec.\ref{subsubsec:comparison_with_QUBO_solvers}}}
    \label{tab:AMFD_parameters}
    % \fontsize{10pt}{0.6cm}\selectfont
\begin{tabular}{ccccc}
\hline
Instance       & $\eta$ & $\zeta$ & $T_\text{init}$ & $T_\text{fin}$ \\ \hline\hline
MCP (G35)      & 0.2    & 5       & 0.3             & 0.1            \\
MISP (C2000.9) & 0.05   & 2       & 0.5             & 0              \\
TSP (eil51)    & 0.05   & 0       & 0.3             & 0.1            \\
QAP (tai50a)   & 0.005  & 0       & 0.3             & 0.1            \\
GCP (jean)     & 0.005  & 50      & 0.3             & 0              \\ \hline
\end{tabular}
\end{table}

In \textbf{Sec.\ref{subsubsec:comparison_with_QUBO_solvers}}, the details of the parameters used for each QUBO solver are shown in \textbf{Tab. \ref{tab:RPA_parameters}--\ref{tab:AMFD_parameters}}.

\textbf{Tab. \ref{tab:RPA_parameters}} shows the parameters of RPA for each COP, where $\varepsilon$ indicates the parameter to prevent excessive flipping, $T_{\text{init}}$ represents the initial temperature, and $T_{\text{fin}}$ represents the final temperature. Following the original paper, we convert the QUBO function represented by 0 and 1 variables into the Ising function represented by -1 and +1 variables, where $J_{i,j}$ represents the interaction in the Ising function and $m_i$ represents the external magnetic field in the Ising function.
Then, the normalization constant $M$, which serves as the base for determining the temperature parameters in the table, is calculated using the following formula:
\begin{equation}\label{eq:RPA_normalization}
    M = \sqrt{\frac{\sum_{i=1}^{N_\text{spin}}(m_i^2 + \sum_{j=1}^{N_\text{spin}}J_{i,j}^2)}{N_\text{spin}}}.
\end{equation}
 $\varepsilon$ is selected from $\{0.1, 0.2, 0.3, 0.4, 0.5, 0.6, 0.7, 0.8, 0.9\}$, $T_{\text{init}}$ is chosen from $\{0.5M, M\}$, and $T_{\text{fin}}$ is chosen from $\{0.5M/N_{\text{spin}}, M/N_{\text{spin}}, 2M/N_{\text{spin}}\}$ to be the best performance.

\textbf{Tab. \ref{tab:bSB_parameters}, \ref{tab:dSB_parameters}} represent the parameters of bSB and dSB, where $\Delta t$ denotes the time width in the time evolution simulation, $a_0$ represents the coefficient of the initial transverse magnetic field, and $c_0$ represents the coefficient applied to the Ising function. Each parameter, based on the original paper, we select $\Delta t$ from $\{0.05, 0.1, 0.25, 0.5, 0.75, 1, 1.25\}$, $a_0$ from $\{0.5, 1\}$, and $c_0$ from $\{0.5/G, 1/G\}$ to obtain the best performance, where $G$ is expressed by
\begin{equation}\label{eq:SB_normalization}
    G = \sqrt{\frac{\sum_{i=1}^{N_\text{spin}}\sum_{j=1}^{N_\text{spin}}J_{i,j}^2}{N_\text{spin} - 1}}.
\end{equation}

\textbf{Tab. \ref{tab:NMFA_parameters}} shows the parameters for NMFA, where $\alpha$ represents the rate of progression towards the values calculated in the self-consistent equation, $\sigma$ is the standard deviation of the Gaussian noise, $T_{\text{init}}$ is the initial temperature, and $T_{\text{fin}}$ is the final temperature. We tune each parameter by selecting $\alpha$ from $\{0.1, 0.2, 0.3, 0.4, 0.5, 0.6, 0.7, 0.8, 0.9\}$, $T_{\text{init}}$ from $\{0.005, 0.01, 0.05, 0.1, 0.5, 1\}$, $T_{\text{fin}}$ from $\{0, 0.0001, 0.001\}$, and $\sigma$ from $\{0.0001, 0.001, 0.01, 0.1\}$.

\textbf{Tab. \ref{tab:VMFA_parameters}} represents the parameters for VMFA, where $\eta$ represents the step size, $T_{\text{init}}$ represents the initial temperature, and $T_{\text{fin}}$ represents the final temperature. Each parameter is selected as follows: $\eta$ from $\{0.05, 0.1, 0.2, 0.5, 0.75, 1\}$, $T_{\text{init}}$ from $\{V, 0.5V\}$, and $T_{\text{fin}}$ from $\{0.5V/N_\text{spin}, V/N_\text{spin}\}$, where 
\begin{equation}
    \label{VMFA_normalization}
    V = \sqrt{\frac{\sum_{i=1}^{N_\text{spin}}(h_i^2 + \sum_{j=1}^{N_\text{spin}}Q_{i,j}^2)}{N_\text{spin}}}
\end{equation}

\textbf{Tab. \ref{tab:AMFD_parameters}} shows the parameters of AMFD, where $\eta$ represents the step size, $\zeta$ denotes the degree of advancement, $T_{\text{init}}$ represents the initial temperature, and $T_{\text{fin}}$ represents the final temperature. Each parameter is chosen as follows: $\eta$ from $\{0.002, 0.005, 0.01, 0.02, 0.05, 0.1, 0.2\}$, $\zeta$ from $\{0,1,2,5,10,20,50\}$, $T_{\text{init}}$ from $\{0.3, 0.5\}$, and $T_{\text{fin}}$ from $\{0, 0.1\}$, optimized for each problem.

\section{Computational Cost Analysis}\label{sec:computational_cost_analysis}
\begin{figure}[b]
    \centering
    \includegraphics[width=0.9\linewidth]{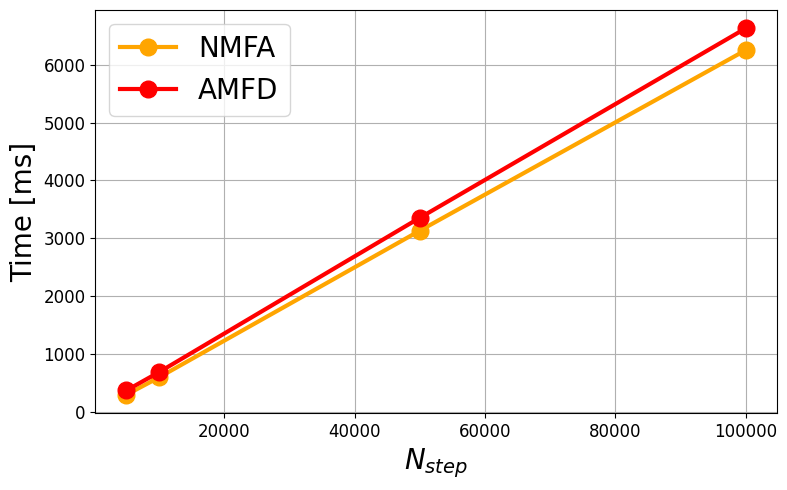}
    \caption{Computation time as a function of $N_\text{step}$}
    \label{fig:computational_time}
\end{figure}

{As discussed in \textbf{Sec.~\ref{subsec:experimental_setting}}, the dominant cost per step in all parallel SAs lies in computing the local or mean field $\boldsymbol{\Phi}$. 
This cost is nearly identical to that of evaluating the gradient of the KL divergence in AMFD. 
Therefore, the asymptotic per-step computational complexity is the same for all algorithms. 
However, the overhead from other processing steps slightly differs among the algorithms, making it necessary to examine their actual contribution to the total computation time. 
Here, we compare the runtime with NMFA, whose computation outside the local-field update is the simplest. 
\textbf{Fig.~\ref{fig:computational_time}} presents the computation times for NMFA and AMFD when $N_\text{step}$ is varied over $\{5000,\, 10000,\, 50000,\, 100000\}$. 
The experiments were conducted on the \textsc{G34} MCP instance with 2000 nodes, running 128 parallel trials on the same GPU, and the total wall-clock time was measured. 
The results indicate that the difference in computation time is marginal, and no substantial increase in runtime is observed.}

\section{Hyperparameter Sensitivity}\label{sec:hyperparameter_sensitivity}

\begin{table}[h]
    \centering
    \footnotesize
    \caption{Best gap [\%] over 100 runs for MCP (\texttt{G22}) under varying $\eta$ and $\zeta$.}
    \label{tab:MCP_sensitivity}
\begin{tabular}{c|rrrrrrr}
\hline
\multicolumn{1}{l|}{$\eta$\textbackslash{}$\zeta$} & \multicolumn{1}{c}{0} & \multicolumn{1}{c}{1} & \multicolumn{1}{c}{2} & \multicolumn{1}{c}{5} & \multicolumn{1}{c}{10} & \multicolumn{1}{c}{20} & \multicolumn{1}{c}{50} \\ \hline
0.002                                        & 0.5                   & 0.5                   & 0.6                   & 0.4                   & 0.1                    & 0.1                    & 0.0                    \\
0.005                                        & 0.4                   & 0.5                   & 0.2                   & 0.3                   & 0.2                    & 0.0                    & 3.1                    \\
0.01                                         & 0.4                   & 0.5                   & 0.2                   & 0.3                   & 0.0                    & 0.0                    & 3.4                    \\
0.02                                         & 0.4                   & 0.3                   & 0.3                   & 0.0                   & 0.0                    & 0.5                    & 4.0                    \\
0.05                                         & 0.0                   & 0.3                   & 0.2                   & 0.0                   & 0.0                    & 3.2                    & 24.3                   \\
0.1                                          & 0.5                   & 0.1                   & 0.0                   & 0.0                   & 0.7                    & 5.0                    & 21.9                   \\
0.2                                          & 0.3                   & 0.0                   & 0.0                   & 0.0                   & 4.4                    & 22.8                   & 28.7                   \\ \hline
\end{tabular}
\end{table}

\begin{table}[h]
\footnotesize
    \centering
    \caption{Best gap [\%] over 100 runs for MISP (\texttt{C2000.5}) under varying $\eta$ and $\zeta$.}
    \label{tab:MISP_sensitivity}
\begin{tabular}{c|rrrrrrr}
\hline
\multicolumn{1}{l|}{$\eta$\textbackslash{}$\zeta$} & \multicolumn{1}{c}{0} & \multicolumn{1}{c}{1} & \multicolumn{1}{c}{2} & \multicolumn{1}{c}{5} & \multicolumn{1}{c}{10} & \multicolumn{1}{c}{20} & \multicolumn{1}{c}{50} \\ \hline
0.002                                        & 0.0                   & 0.0                   & 6.3                   & 0.0                   & 0.0                    & 0.0                    & 6.3                    \\
0.005                                        & 6.3                   & 6.3                   & 0.0                   & 0.0                   & 0.0                    & 6.3                    & 1.0e7             \\
0.01                                         & 6.3                   & 6.3                   & 6.3                   & 0.0                   & 0.0                    & 12.5                   & 1.2e7            \\
0.02                                         & 0.0                   & 0.0                   & 0.0                   & 0.0                   & 1.2e7            & 1.2e7            & 1.2e7            \\
0.05                                         & 6.3                   & 0.0                   & 18.8                  & 12.5                  & 18.8                   & 31.3                   & 50.0                   \\
0.1                                          & 12.5                  & 37.5                  & 18.8                  & 62.5                  & 100.0                  & 100.0                  & 100.0                  \\
0.2                                          & 31.3                  & 25.0                  & 43.8                  & 100.0                 & 100.0                  & 100.0                  & 100.0                  \\ \hline
\end{tabular}
\end{table}

\begin{table}[h]
\footnotesize
\centering
\caption{Best gap [\%] over 100 runs for TSP (\texttt{eil51}) under varying $\eta$ and $\zeta$.}
\label{tab:tsp_sensitivity}
\begin{tabular}{c|rrrrrrr}
\hline
\multicolumn{1}{l|}{$\eta$\textbackslash{}$\zeta$} & \multicolumn{1}{c}{0} & \multicolumn{1}{c}{1} & \multicolumn{1}{c}{2} & \multicolumn{1}{c}{5} & \multicolumn{1}{c}{10} & \multicolumn{1}{c}{20} & \multicolumn{1}{c}{50} \\ \hline
0.002                                        & 3.1                   & 4.2                   & 3.1                   & 5.9                   & 7.7                    & 9.4                    & 17.4                   \\
0.005                                        & 2.1                   & 3.3                   & 3.8                   & 5.2                   & 7.5                    & 13.6                   & 29.8                   \\
0.01                                         & 3.5                   & 2.8                   & 3.3                   & 3.8                   & 10.1                   & 15.7                   & 119.3                  \\
0.02                                         & 3.1                   & 4.0                   & 3.5                   & 6.8                   & 11.3                   & 54.0                   & 3.7e6              \\
0.05                                         & 0.9                   & 2.8                   & 2.8                   & 18.1                  & 3.7e6              & 3.7e6              & 3.7e6              \\
0.1                                          & 7.7                   & 5.2                   & 20.9                  & 70.0                  & 2.8e4                & 8.9e3                 & 5.6e3                 \\
0.2                                          & 582.2                 & 599.3                 & 3.2e4               & 5.5e5              & 1.4e6              & 1.0e3                 & 4.3e3                 \\ \hline
\end{tabular}
\end{table}

\begin{table}[h]
\footnotesize
\centering
\caption{Best gap [\%] over 100 runs for QAP (\texttt{tai50a}) under varying $\eta$ and $\zeta$.}
\label{tab:qap_sensitivity}
\begin{tabular}{c|rrrrrrr}
\hline
\multicolumn{1}{l|}{$\eta$\textbackslash{}$\zeta$} & \multicolumn{1}{c}{0} & \multicolumn{1}{c}{1} & \multicolumn{1}{c}{2} & \multicolumn{1}{c}{5} & \multicolumn{1}{c}{10} & \multicolumn{1}{c}{20} & \multicolumn{1}{c}{50} \\ \hline
0.002                                        & 2.4                   & 2.5                   & 2.5                   & 2.5                   & 3.0                    & 6.3                    & 7.5                    \\
0.005                                        & 2.4                   & 1.8                   & 2.2                   & 2.3                   & 6.2                    & 7.2                    & 9.4                    \\
0.01                                         & 3.5                   & 3.4                   & 3.8                   & 4.5                   & 7.0                    & 8.3                    & 10.5                   \\
0.02                                         & 5.1                   & 5.2                   & 5.1                   & 6.9                   & 8.6                    & 9.5                    & 1.1e6             \\
0.05                                         & 7.2                   & 7.1                   & 6.5                   & 10.4                  & 1.1e6             & 1.1e6             & 1.1e6             \\
0.1                                          & 8.0                   & 7.6                   & 9.6                   & 9.3                   & 10.5                   & 225.8                  & 944.2                  \\
0.2                                          & 10.7                  & 10.3                  & 10.1                  & 232.4                 & 232.4                  & 232.4                  & 796.7                  \\ \hline
\end{tabular}
\end{table}

\begin{table}[h]
\caption{Best gap [\%] over 100 runs for GCP (\texttt{jean}) under varying $\eta$ and $\zeta$.}
\label{tab:gcp_sensitivity}
\centering
\begin{tabular}{c|rrrrrrr}
\hline
\multicolumn{1}{l|}{$\eta$\textbackslash{}$\zeta$} & \multicolumn{1}{c}{0} & \multicolumn{1}{c}{1} & \multicolumn{1}{c}{2} & \multicolumn{1}{c}{5} & \multicolumn{1}{c}{10} & \multicolumn{1}{c}{20} & \multicolumn{1}{c}{50} \\ \hline
0.002                                        & 100                   & 100                   & 90                    & 80                    & 80                     & 70                     & 10                     \\
0.005                                        & 80                    & 80                    & 70                    & 80                    & 60                     & 40                     & 0                      \\
0.01                                         & 70                    & 70                    & 70                    & 70                    & 60                     & 20                     & 1.3e4                 \\
0.02                                         & 50                    & 50                    & 60                    & 60                    & 40                     & 0                      & 1.1e6                \\
0.05                                         & 20                    & 40                    & 50                    & 40                    & 30                     & 1.5e6                & 2.3e6                \\
0.1                                          & 30                    & 40                    & 50                    & 70                    & 2.3e6                & 2.3e6                & 2.8e5                 \\
0.2                                          & 60                    & 80                    & 70                    & 120                   & 2.4e5                 & 3.2e5                 & 4.0e5                \\ \hline
\end{tabular}
\end{table}

{The AMFD algorithm involves four hyperparameters: the step size $\eta$, the inertial advancement $\zeta$, the initial temperature $T_\text{init}$, and the final temperature $T_\text{fin}$. Among these, $T_\text{init}$ and $T_\text{fin}$ have limited impact on solution quality when the number of iterations $N_\text{step}$ is sufficiently large, as they do not induce abrupt temperature changes. In contrast, $\eta$ and $\zeta$ have a pronounced influence on performance, with solution quality varying significantly across parameter choices. 

\textbf{Tab.~\ref{tab:MCP_sensitivity}, \ref{tab:MISP_sensitivity}, \ref{tab:tsp_sensitivity}, \ref{tab:qap_sensitivity}, \ref{tab:gcp_sensitivity}} report the best gap [\%] obtained over 100 independent runs for each grid point, where $\eta \in \{0.002, 0.005, 0.01, 0.02, 0.05, 0.1, 0.2\}$ and $\zeta \in \{0, 1, 2, 5, 10, 20, 50\}$. The evaluation covers MCP (\texttt{G22}), MISP (\texttt{C2000.5}), TSP (\texttt{eil51}), QAP (\texttt{tai50a}), and GCP (\texttt{jean}). 
The results indicate that, once $\eta$ or $\zeta$ exceed certain thresholds, convergence deteriorates sharply, leading to degraded solution quality. Below these thresholds, a broad range of parameter settings still yield reasonably good solutions. Notably, performance tends to improve when the product $\eta \cdot \zeta$ is pushed as close as possible to the non-convergent boundary. This suggests that, although the exact Pareto frontier varies across problem instances, identifying it is crucial for achieving high-quality solutions. 

As demonstrated in this paper, exploring multiple parameter configurations in parallel is computationally feasible and does not impose a major bottleneck. Nevertheless, in future work, leveraging machine learning or statistical inference to promptly estimate near-optimal parameters could be valuable for reducing runtime while preserving solution quality.}

\section{Detailed Experimental Results}\label{sec:detailed_results}

{
This section presents the detailed results of the experiments described in \textbf{Sec.\ref{subsec:experimental_results}}. In \textbf{Tab. \ref{tab:detailed_MCP}}, \textbf{\ref{tab:detailed_MISP}}, \textbf{\ref{tab:detailed_TSP}}, \textbf{\ref{tab:detailed_QAP}}, and \textbf{\ref{tab:detailed_GCP}}, for each instance, the first row in the AMFD column reports the solution and execution time obtained in the tuning phase. The subsequent five rows correspond to the refinement phase, where $N_\text{step}$ was set to $2\times$, $5\times$, $10\times$, $20\times$, and $50\times N_\text{spin}$, respectively; for each setting, the table reports the total execution time (including the tuning phase) and the solution obtained. For Gurobi, the time limit was set to twice the maximum execution time of AMFD (i.e., twice the sum of the tuning phase time and the time for $N_\text{step} = 50N_\text{spin}$ in the refinement phase), and the table lists representative solutions obtained within this time limit along with the times at which they were first found.
}

\clearpage
\onecolumn
\begin{center}
\centering
\begin{longtable}{crrrrrrrrr}
\caption{Detailed results of experiments on MCP instances.} \label{tab:detailed_MCP} 
\\
\hline
\multirow{2}{*}{Instance} & \multicolumn{1}{c}{\multirow{2}{*}{BKS}} 
& \multicolumn{2}{c}{Gurobi (MILP)} 
& \multicolumn{2}{c}{Gurobi (MIQP)}  
& \multicolumn{2}{c}{Gurobi (QUBO)}  
& \multicolumn{2}{c}{AMFD (QUBO)} \\

& \multicolumn{1}{c}{}                     
& \multicolumn{1}{c}{Time [s]} & \multicolumn{1}{c}{SOL} 
& \multicolumn{1}{c}{Time [s]} & \multicolumn{1}{c}{SOL} 
& \multicolumn{1}{c}{Time [s]} & \multicolumn{1}{c}{SOL} 
& \multicolumn{1}{c}{Time [s]} & \multicolumn{1}{c}{SOL} \\
\hline\hline

\endfirsthead

\hline
\multirow{2}{*}{Instance} & \multicolumn{1}{c}{\multirow{2}{*}{BKS}} 
& \multicolumn{2}{c}{Gurobi (MILP)} 
& \multicolumn{2}{c}{Gurobi (MIQP)}  
& \multicolumn{2}{c}{Gurobi (QUBO)}  
& \multicolumn{2}{c}{AMFD (QUBO)} \\

& \multicolumn{1}{c}{}                     
& \multicolumn{1}{c}{Time [s]} & \multicolumn{1}{c}{SOL} 
& \multicolumn{1}{c}{Time [s]} & \multicolumn{1}{c}{SOL} 
& \multicolumn{1}{c}{Time [s]} & \multicolumn{1}{c}{SOL} 
& \multicolumn{1}{c}{Time [s]} & \multicolumn{1}{c}{SOL} \\
\hline\hline
\endhead

G1                        & 11624                                    & 0.009                    & 0                            &   N/A                       &  N/A                            & 0.001                    & 0                            & 0.853                    & 11623                        \\
G1                        & 11624                                    & 0.317                    & 9426                         &                          &                              & 0.033                    & 11514                        & 2.437                    & 11624                        \\
G1                        & 11624                                    & 0.766                    & 10175                        &                          &                              & 0.041                    & 11593                        & 2.826                    & 11624                        \\
G1                        & 11624                                    & 1.941                    & 10278                        &                          &                              & 8.871                    & 11594                        & 3.230                    & 11624                        \\
G1                        & 11624                                    & 3.496                    & 10297                        &                          &                              & 13.690                   & 11597                        & 3.993                    & 11624                        \\
G1                        & 11624                                    & 3.633                    & 10306                        &                          &                              & 14.091                   & 11598                        & 7.772                    & 11624                        \\
\hline

\end{longtable}
\end{center}

\begin{center}
\begin{longtable}{crrrrrrrrr}
\caption{Detailed results of experiments on MISP instances.} \label{tab:detailed_MISP} \\

\hline
\multirow{2}{*}{Instance} & \multicolumn{1}{c}{\multirow{2}{*}{BKS}} 
& \multicolumn{2}{c}{Gurobi (MILP)} 
& \multicolumn{2}{c}{Gurobi (MIQP)}  
& \multicolumn{2}{c}{Gurobi (QUBO)}  
& \multicolumn{2}{c}{AMFD (QUBO)} \\

& \multicolumn{1}{c}{}                     
& \multicolumn{1}{c}{Time [s]} & \multicolumn{1}{c}{SOL} 
& \multicolumn{1}{c}{Time [s]} & \multicolumn{1}{c}{SOL} 
& \multicolumn{1}{c}{Time [s]} & \multicolumn{1}{c}{SOL} 
& \multicolumn{1}{c}{Time [s]} & \multicolumn{1}{c}{SOL} \\
\hline\hline

\endfirsthead

\hline
\multirow{2}{*}{Instance} & \multicolumn{1}{c}{\multirow{2}{*}{BKS}} 
& \multicolumn{2}{c}{Gurobi (MILP)} 
& \multicolumn{2}{c}{Gurobi (MIQP)}  
& \multicolumn{2}{c}{Gurobi (QUBO)}  
& \multicolumn{2}{c}{AMFD (QUBO)} \\

& \multicolumn{1}{c}{}                     
& \multicolumn{1}{c}{Time [s]} & \multicolumn{1}{c}{SOL} 
& \multicolumn{1}{c}{Time [s]} & \multicolumn{1}{c}{SOL} 
& \multicolumn{1}{c}{Time [s]} & \multicolumn{1}{c}{SOL} 
& \multicolumn{1}{c}{Time [s]} & \multicolumn{1}{c}{SOL} \\
\hline\hline
\endhead

brock200\_1               & 21                                       & 0.0003                   & 13                           & 0.001                    & 0                            & 0.000                    & 0                            & 1.482                    & 20                           \\
brock200\_1               & 21                                       & 0.173                    & 15                           & 0.059                    & 13                           & 0.028                    & 19                           & 3.059                    & 21                           \\
brock200\_1               & 21                                       & 0.351                    & 17                           & 0.203                    & 18                           & 0.660                    & 20                           & 3.154                    & 21                           \\
brock200\_1               & 21                                       & 2.345                    & 19                           & 0.556                    & 19                           & 1.559                    & 21                           & 3.227                    & 21                           \\
brock200\_1               & 21                                       & 4.673                    & 21                           & 4.374                    & 21                           & 1.621                    & 21                           & 3.588                    & 21                           \\
brock200\_1               & 21                                       & 4.908                    & 21                           & 4.577                    & 21                           &                          &                              & 4.421                    & 21                           \\
\hline

\end{longtable}
\end{center}

\begin{center}
\begin{longtable}{crrrrrrrrr}
\caption{Detailed results of experiments on TSP instances.} \label{tab:detailed_TSP} \\

\hline
\multirow{2}{*}{Instance} & \multicolumn{1}{c}{\multirow{2}{*}{BKS}} 
& \multicolumn{2}{c}{Gurobi (MILP)} 
& \multicolumn{2}{c}{Gurobi (MIQP)}  
& \multicolumn{2}{c}{Gurobi (QUBO)}  
& \multicolumn{2}{c}{AMFD (QUBO)} \\

& \multicolumn{1}{c}{}                     
& \multicolumn{1}{c}{Time [s]} & \multicolumn{1}{c}{SOL} 
& \multicolumn{1}{c}{Time [s]} & \multicolumn{1}{c}{SOL} 
& \multicolumn{1}{c}{Time [s]} & \multicolumn{1}{c}{SOL} 
& \multicolumn{1}{c}{Time [s]} & \multicolumn{1}{c}{SOL} \\
\hline\hline

\endfirsthead

\hline
\multirow{2}{*}{Instance} & \multicolumn{1}{c}{\multirow{2}{*}{BKS}} 
& \multicolumn{2}{c}{Gurobi (MILP)} 
& \multicolumn{2}{c}{Gurobi (MIQP)}  
& \multicolumn{2}{c}{Gurobi (QUBO)}  
& \multicolumn{2}{c}{AMFD (QUBO)} \\

& \multicolumn{1}{c}{}                     
& \multicolumn{1}{c}{Time [s]} & \multicolumn{1}{c}{SOL} 
& \multicolumn{1}{c}{Time [s]} & \multicolumn{1}{c}{SOL} 
& \multicolumn{1}{c}{Time [s]} & \multicolumn{1}{c}{SOL} 
& \multicolumn{1}{c}{Time [s]} & \multicolumn{1}{c}{SOL} \\
\hline\hline
\endhead

% Please add the following required packages to your document preamble:
% \usepackage{multirow}

att48                     & 10628                                    & 0.067                    & 67822                        & 0.003                    & 49358                        & 0.837                    & 18013                        & 1.544                    & 11278                        \\
att48                     & 10628                                    & 0.268                    & 12156                        & 13.650                   & 31051                        & 0.934                    & 16731                        & 3.478                    & 10858                        \\
att48                     & 10628                                    & 1.584                    & 10821                        & 17.729                   & 17338                        & 9.907                    & 16730                        & 4.000                    & 10739                        \\
att48                     & 10628                                    & 2.043                    & 10704                        & 26.939                   & 17296                        & 40.984                   & 16730                        & 5.788                    & 10805                        \\
att48                     & 10628                                    & 4.137                    & 10653                        & 31.627                   & 16418                        &                          &                              & 9.310                    & 10729                        \\
att48                     & 10628                                    & 4.278                    & 10628                        & 35.120                   & 16104                        &                          &                              & 20.145                   & 10729                        \\
\hline

\end{longtable}
\end{center}

\begin{center}
\begin{longtable}{crrrrrrrrr}
\caption{Detailed results of experiments on QAP instances.} \label{tab:detailed_QAP} \\

\hline
\multirow{2}{*}{Instance} & \multicolumn{1}{c}{\multirow{2}{*}{BKS}} 
& \multicolumn{2}{c}{Gurobi (MILP)} 
& \multicolumn{2}{c}{Gurobi (MIQP)}  
& \multicolumn{2}{c}{Gurobi (QUBO)}  
& \multicolumn{2}{c}{AMFD (QUBO)} \\

& \multicolumn{1}{c}{}                     
& \multicolumn{1}{c}{Time [s]} & \multicolumn{1}{c}{SOL} 
& \multicolumn{1}{c}{Time [s]} & \multicolumn{1}{c}{SOL} 
& \multicolumn{1}{c}{Time [s]} & \multicolumn{1}{c}{SOL} 
& \multicolumn{1}{c}{Time [s]} & \multicolumn{1}{c}{SOL} \\
\hline\hline

\endfirsthead

\hline
\multirow{2}{*}{Instance} & \multicolumn{1}{c}{\multirow{2}{*}{BKS}} 
& \multicolumn{2}{c}{Gurobi (MILP)} 
& \multicolumn{2}{c}{Gurobi (MIQP)}  
& \multicolumn{2}{c}{Gurobi (QUBO)}  
& \multicolumn{2}{c}{AMFD (QUBO)} \\

& \multicolumn{1}{c}{}                     
& \multicolumn{1}{c}{Time [s]} & \multicolumn{1}{c}{SOL} 
& \multicolumn{1}{c}{Time [s]} & \multicolumn{1}{c}{SOL} 
& \multicolumn{1}{c}{Time [s]} & \multicolumn{1}{c}{SOL} 
& \multicolumn{1}{c}{Time [s]} & \multicolumn{1}{c}{SOL} \\
\hline\hline
\endhead

bur26a                    & 5426670                                  & 0.007                    & 6147379                      & 0.002                    & 6147379                      & 0.798                    & 5599040                      & 0.992                    & 5644287                      \\
bur26a                    & 5426670                                  & 1.303                    & 6053533                      & 4.336                    & 5728600                      & 2.585                    & 5556103                      & 2.305                    & 5497748                      \\
bur26a                    & 5426670                                  & 2.888                    & 5724504                      & 4.466                    & 5648700                      & 8.904                    & 5533873                      & 2.439                    & 5545402                      \\
bur26a                    & 5426670                                  & 4.302                    & 5674468                      & 10.182                   & 5533878                      & 9.848                    & 5529343                      & 2.778                    & 5505928                      \\
bur26a                    & 5426670                                  & 7.226                    & 5646231                      & 12.117                   & 5528569                      & 9.906                    & 5514267                      & 3.863                    & 5484383                      \\
bur26a                    & 5426670                                  & 13.276                   & 5621263                      & 14.129                   & 5528569                      & 9.959                    & 5500764                      & 7.025                    & 5496503                      \\
\hline

\end{longtable}
\end{center}

\begin{center}
\begin{longtable}{crrrrrrrrr}
\caption{Detailed results of experiments on GCP instances.} \label{tab:detailed_GCP} \\

\hline
\multirow{2}{*}{Instance} & \multicolumn{1}{c}{\multirow{2}{*}{BKS}} 
& \multicolumn{2}{c}{Gurobi (MILP)} 
& \multicolumn{2}{c}{Gurobi (MIQP)}  
& \multicolumn{2}{c}{Gurobi (QUBO)}  
& \multicolumn{2}{c}{AMFD (QUBO)} \\

& \multicolumn{1}{c}{}                     
& \multicolumn{1}{c}{Time [s]} & \multicolumn{1}{c}{SOL} 
& \multicolumn{1}{c}{Time [s]} & \multicolumn{1}{c}{SOL} 
& \multicolumn{1}{c}{Time [s]} & \multicolumn{1}{c}{SOL} 
& \multicolumn{1}{c}{Time [s]} & \multicolumn{1}{c}{SOL} \\
\hline\hline

\endfirsthead

\hline
\multirow{2}{*}{Instance} & \multicolumn{1}{c}{\multirow{2}{*}{BKS}} 
& \multicolumn{2}{c}{Gurobi (MILP)} 
& \multicolumn{2}{c}{Gurobi (MIQP)}  
& \multicolumn{2}{c}{Gurobi (QUBO)}  
& \multicolumn{2}{c}{AMFD (QUBO)} \\

& \multicolumn{1}{c}{}                     
& \multicolumn{1}{c}{Time [s]} & \multicolumn{1}{c}{SOL} 
& \multicolumn{1}{c}{Time [s]} & \multicolumn{1}{c}{SOL} 
& \multicolumn{1}{c}{Time [s]} & \multicolumn{1}{c}{SOL} 
& \multicolumn{1}{c}{Time [s]} & \multicolumn{1}{c}{SOL} \\
\hline\hline
\endhead
                          
anna                      & 11                                       & 0.008                    & 62                           & 0.423                    & 72                           & 1.612                    & 15                           & 8.152                    & 11                           \\
anna                      & 11                                       & 0.272                    & 61                           & 0.432                    & 71                           & 1.845                    & 12                           & 11.679                   & 11                           \\
anna                      & 11                                       & 3.229                    & 11                           & 0.445                    & 12                           & 55.187                   & 12                           & 12.547                   & 11                           \\
anna                      & 11                                       & 3.239                    & 11                           & 1.158                    & 11                           &                          &                              & 14.324                   & 11                           \\
anna                      & 11                                       &                          &                              & 1.302                    & 11                           &                          &                              & 17.409                   & 11                           \\
anna                      & 11                                       &                          &                              &                          &                              &                          &                              & 26.906                   & 11                           \\
\hline

\end{longtable}
\end{center}
% \clearpage
% \pagestyle{empty} 
% \input{9_figures_and_tables}

\end{document}